\documentclass[a4paper]{article}

\usepackage{amsfonts, amsmath, amsthm, empheq, amssymb, color, indentfirst, physics}

\allowdisplaybreaks[4]
\newtheorem{theorem}{Theorem}[section]
\newtheorem{proposition}{Proposition}[section]
\newtheorem{corollary}{Corollary}[section]
\newtheorem{lemma}{Lemma}[section]
\newtheorem{remark}{Remark}[section]

\newcommand{\OOO}{\Omega}
\newcommand{\ooo}{\overline}
\newcommand{\ppp}{\partial}

\newcommand{\LLLLLL}{L^2(0,T;X)}
\newcommand{\Halp}{H_{\alpha}(0,T;X)}
\newcommand{\Hgam}{H_{\gamma}(0,T;X)}
\newcommand{\pppa}{\partial_t^{\alpha}}
\newcommand{\pppg}{\partial_t^{\gamma}}
\newcommand{\DDD}{\mathcal{D}}
\newcommand{\N}{\mathbb{N}}
\newcommand{\R}{\mathbb{R}}
\newcommand{\JJJa}{J^{\gamma}}
\newcommand{\ep}{\varepsilon}
\newcommand{\pppppp}{\frac{\ppp}{\ppp t}}
\newcommand{\dddddd}{\frac{d}{dt}}
\newcommand{\pun}{\frac{\ppp u_N}{\ppp t}}

\title{\Large\bf
Well-posedness of initial-boundary value problem for time-fractional diffusion-wave equation 
\\with time-dependent coefficients}

\author{\large Xinchi HUANG$^{1,\ast}$, Masahiro YAMAMOTO$^2$}

\date{}

\begin{document}
\maketitle

\renewcommand{\thefootnote}{\fnsymbol{footnote}}
\footnotetext{\hspace*{-5mm} 
\begin{tabular}{@{}r@{}p{14cm}@{}} 
& Manuscript last updated: March 12, 2025.\\
$^1$ 
& School of Science, The University of Tokyo,
7-3-1 Hongo, Bunkyo-ku, Tokyo 113-8654, Japan.
E-mail: 
huangxc@g.ecc.u-tokyo.ac.jp\\
$^2$ 
& Graduate School of Mathematical Sciences, 
the University of Tokyo,
3-8-1 Komaba, Meguro-ku, Tokyo 153-8914, Japan. \\
& Department of Mathematics, Faculty of Science, Zonguldak B\"ulent Ecevit University, 
Zonguldak 67100, T\"urkiye.\\
& Honorary Member of Academy of Romanian Scientists,
Ilfov, nr. 3, Bucuresti, Romania.\\
& Correspondence member of Accademia Peloritana dei Pericolanti, Palazzo Universit\`a, 
Piazza S. Pugliatti 1 98122 Messina, Italy. 
E-mail: myama@ms.u-tokyo.ac.jp\\
$^\ast$
& Corresponding author: huangxc@g.ecc.u-tokyo.ac.jp
\end{tabular}}

\begin{abstract}
We consider the well-posedness of the initial-boundary value problem for a time-fractional partial differential equation with the fractional order lying in (1,2]. 
For the case of time-dependent coefficients, it is difficult to give an explicit solution formula by the eigenfunction expansion method. 
In order to deal with the case of time-varying coefficients, we first show the unique existence and regularity of solution to a system of time-fractional ordinary differential equations. 
Then the unique existence of the weak solution to the time-fractional partial differential equation and improved regularity are derived by using the Galerkin method.

\vskip 4.5mm

\noindent\begin{tabular}{@{}l@{ }p{11.5cm}} {\bf Keywords } &
Time-fractional diffusion/wave equation,
Initial-boundary value problem, 
Fredholm alternative,
Galerkin approximation,
Regularity estimate
\end{tabular}

\vskip 4.5mm


\end{abstract}

\section{Introduction and main results}
\label{sec-intro}
\subsection{Introduction}

During the recent decades, there have been incredible growing interests on the topic of the time-fractional diffusion/wave equation: 
$$
\partial_t^\alpha u + A u = F, \quad 0<\alpha\le 2,
$$
where $A$ is a second-order elliptic operator, and $\partial_t^\alpha$ denotes the fractional derivative 
which will be introduced later. According to the background of its application, we may call it a time-fractional diffusion equation while $0<\alpha<1$, and a time-fractional diffusion-wave equation while $1<\alpha<2$. 

From the mathematical viewpoint, the time-fractional diffusion/wave equation is a 
generalization of the classical diffusion equation or the wave equation since one considers 
an integro-differential operator of real order $\alpha\in (0,2]$ and its solution corresponds to 
the solutions of the diffusion equation as $\alpha\to 1$ and the wave equation 
as $\alpha\to 2$, respectively. 
From the practical viewpoint, the time-fractional diffusion/wave equation is regarded as 
a feasible candidate for modeling the anomalous diffusion in porous media 
(see e.g., Nigmatullin \cite{N86} for $0<\alpha<1$) or the mechanical diffusive waves in 
viscoelastic media which admit a power-law creep (see e.g., Mainardi \cite{M95} for $1<\alpha<2$). 
Actually with the time-fractional derivatives, one can expect some memory effects of power-law type 
in the diffusion or the wave phenomena. 

Owing to the applications in engineering and other applied sciences, 
the initial-boundary value problems for the time-fractional diffusion/wave 
equations have been intensively investigated in the last few decades. 
For the special cases of time-independent coefficients, 
there are many existing results. 
For example, one can readily derive the solution formula to the time-fractional diffusion/wave equation 
in terms of the Fourier method and the so-called Mittag-Leffler functions 
\cite{SY11}. 
For some generalizations of the equation including the multi-term 
and the distributed-order time-fractional 
diffusion equations with time-independent coefficients, we refer to e.g., 
\cite{K08,LLY15,L11,LY16,MMP08}.  In particular, Theorem 2.7 
in \cite{KY} establishes the unique existence of strong solutions 
for $1<\alpha<2$ in 
different Sobolev spaces. 

On the other hand, there are very few works on the time-fractional 
diffusion/wave equations with time-dependent coefficients,
where the available methods are limited.  
For example, one can no longer apply the Laplace transform and the eigenfunction expansion to obtain an explicit representation of the solution. 
If the principal coefficients $a_{ij}$ in the elliptic part $A$ do not depend on the time variable, 
then one can regard the lower order terms with time-dependent coefficients as a new source term and apply the fixed-point theorem to overcome the difficulty (e.g., \cite{GLY15}). 
However, this is impossible for the general $a_{ij}=a_{ij}(x,t)$. 
In order to deal with the case of time-dependent coefficients, 
Zacher \cite{Z09}, and Kubica, Ryszewska and Yamamoto \cite{KRY20} applied 
the Galerkin approximation method to prove the unique existence of the solution to 
the initial-boundary value problem for the time-fractional diffusion equation ($0<\alpha<1$). 
 
Comparing to the time-fractional diffusion equation, to the authors' best knowledge, 
there are even fewer publications for the time-fractional diffusion-wave equations for $1<\alpha<2$ 
with time-dependent coefficients. 
Here we should mention the articles \cite{E96,Z06},  
in which the authors considered strong solutions to a related parabolic integro-differential equation:
$$
\frac{\partial}{\partial t} u(x,t) + \int_0^t (t-s)^{-\beta} Au(x,s)ds = f, \quad 0<\beta<1,
$$
under some restrictions on $\beta$ and the nonhomogeneous term $f$. 
Roughly speaking, the above equation changes to a time-fractional diffusion-wave equation 
if one takes a time-fractional derivative on both sides.
In particular, \cite{Z06} considers non-homogeneous boundary
conditions (e.g., Theorem 3.4). 
As for other recent works, Han, Kim and Park \cite{HKP20,P23} 
considered the unique solvability and regularity 
in Sobolev spaces with some special weights provided that the initial condition is vanishing. 

The main purpose of this article is to establish the fundamental theory of the well-posedness of 
the initial-boundary value problem for the time-fractional diffusion-wave equation 
(i.e., $1<\alpha<2$) with time-dependent coefficients. 
In addition, we will give some specific regularity estimates including the improved regularity 
in some fractional Sobolev space, 
which also shows a connection with the classical regularity results for parabolic and hyperbolic equations.  
Before the statement, we formulate the fractional derivatives in Section \ref{subsec:fractional-space}. 

\subsection{Time-fractional derivatives in function spaces}
\label{subsec:fractional-space}

For the statement of our main results, following \cite{KRY20}, 
we introduce notations, operators, and function spaces. 
Here and henceforth, by $L^p(D)$ and $W^{k,p}(D)$ with $k\in \mathbb{N}$
and $1\le p \le \infty$,  
we mean the usual Lebesgue spaces and the $k$-th order Sobolev spaces 
of $L^p$ functions in an interval $D\subset (0,T)$ or $D = \Omega$, 
and in particular, for $p=2$ we denote 
$W^{k,2}(D)$ by $H^k(D)$ (see e.g., \cite{AF75,E98}). 

Henceforth, for a Banach space $X$,
we consider a function $t\in (0,T) \, \longrightarrow\,u(t) \in X$.
This interpretation is common, and see e.g., pp.351--352 in \cite{E98},
\S2 in Chapter 1 of Yagi \cite{Yagi}.
Then, we can naturally define spaces such as $L^2(0,T;X)$, and 
it is well-known that $L^2(0,T;X)$ is also a Banach space provided that 
$X$ is a Banach space.

Moreover, for $\gamma>0$, we define
the Riemann-Liouville fractional integral operator $J^{\gamma}$ 
for an $X$-valued function $u \in L^2(0,T;X)$ by 
$$
J^{\gamma}u(t) 
:= \frac{1}{\Gamma(\gamma)}
\int^t_0 (t-s)^{\gamma-1} u(s) ds, 
$$
and we can regard $J^{\gamma}u$ as a function from $(0,T)$ to $X$.  
We set $J^0 = I$: the identity operator on $X$.

Then we show 
\begin{lemma}
\label{lem:J1}
Let $\gamma>0$. 

\noindent $\mathrm{(i)}$ The operator $J^{\gamma}: L^2(0,T;X) \longrightarrow L^2(0,T;X)$ is 
bounded.

\noindent $\mathrm{(ii)}$ The operator $J^{\gamma}: L^2(0,T;X) \longrightarrow L^2(0,T;X)$ is 
injective.

\noindent $\mathrm{(iii)}$ $J^{\alpha+\beta} = J^{\alpha}J^{\beta}$ in $\LLLLLL$ for 
$\alpha, \beta>0$.
\end{lemma}
\begin{proof}
(i) Young's convolution inequality yields
\begin{align*}
&\quad \Vert J^{\gamma}u\Vert_{\LLLLLL}
= \frac{1}{\Gamma(\gamma)} \left\Vert \int^t_0 (t-s)^{\gamma-1}
u(s) ds\right\Vert_{\LLLLLL}\\
&\le \frac{1}{\Gamma(\gamma)} \left\Vert \int^t_0 (t-s)^{\gamma-1}
\Vert u(s)\Vert_X ds \right\Vert_{L^2(0,T)}
\le \frac{1}{\Gamma(\gamma)} \Vert t^{\gamma-1}\Vert_{L^1(0,T)}
\left( \int^T_0 \Vert u(t)\Vert_{X}^2 dt \right)^{\frac{1}{2}} \\
&\le C\Vert u\Vert_{\LLLLLL}.
\end{align*}

\noindent (ii) 
Let $m$ be the integer satisfying $m\le \gamma<m+1$, 
so that $-\gamma+m\in (-1,0]$. 
We assume that $J^{\gamma}u = 0$ in $\LLLLLL$.
Then for almost all $t\in (0,T)$ and any $\varphi\in X^\prime$:
the dual space of $X$, we have
$$
\frac{1}{\Gamma(m+1-\gamma)}
\int^t_0 (t-s)^{-\gamma+m}{_{X^\prime}\langle \varphi, J^{\gamma} u(s)\rangle_X} ds = 0.
$$
Applying Fubini's theorem to a scalar-valued function 
$(t-s)^{-\gamma+m}(s-\xi)^{\gamma-1}{_{X^\prime}\langle \varphi, u(\xi)\rangle_X}$, 
we can calculate as follows.
\begin{align*}
0 &= \frac{1}{\Gamma(m+1-\gamma)\Gamma(\gamma)}
\left(
\int^t_0 (t-s)^{-\gamma+m}\left( \int^s_0 (s-\xi)^{\gamma-1} {_{X^\prime}\langle \varphi, u(\xi)\rangle_X} d\xi
\right) ds \right)\\
&= \frac{1}{\Gamma(m+1-\gamma)\Gamma(\gamma)}
\int^t_0
\left( \int^t_{\xi} (t-s)^{-\gamma+m}(s-\xi)^{\gamma-1} ds \right)
{_{X^\prime}\langle \varphi, u(\xi)\rangle_X} d\xi \\
&=  \frac{1}{\Gamma(m+1-\gamma)\Gamma(\gamma)}
\int^t_0 \Gamma(\gamma)\Gamma(m+1-\gamma) (t-\xi)^m{_{X^\prime}\langle \varphi, u(\xi)\rangle_X} d\xi \\
&= \int^t_0 (t-\xi)^m {_{X^\prime}\langle \varphi, u(\xi)\rangle_X} d\xi,
\end{align*}
which implies that ${_{X^\prime}\langle \varphi, u(t)\rangle_X}=0$ for almost 
all $t \in (0,T)$ and any $\varphi \in X^\prime$ by repeating the 
differentiation in $t$.  Thus, $u=0$ in $L^2(0,T;X)$. 

We can prove (iii) by calculating $J^{\alpha}J^{\beta}$ and 
exchanging the orders of the integrals.
Thus, the proof of Lemma \ref{lem:J1} is complete. 
\end{proof}

Therefore, the inverse to $J^{\gamma}$ exists and is an operator from 
$J^{\gamma}\LLLLLL$ to $\LLLLLL$.  By $J^{-\gamma}$, we denote the inverse 
operator which is understood algebraically.
Then we can define
\begin{align}
\label{intro:eq2-1}
\left\{ 
\begin{array}{rl}
& H_{\gamma}(0,T;X) := J^{\gamma}\LLLLLL, \\
& \Vert v\Vert_{\Hgam}:= \Vert J^{-\gamma}v\Vert_{\LLLLLL}
\quad \mbox{for } \ v\in H_\gamma(0,T;X).
\end{array}
\right.
\end{align}
In other words, we have $\Vert J^{\gamma} u\Vert_{\Hgam}:= \Vert u\Vert_{\LLLLLL}$ for $u\in L^2(0,T;X)$. 
By \eqref{intro:eq2-1}, we easily see that
$$
\Halp \subset H_{\beta}(0,T;X) \,\, \mbox{algebraically 
if $0<\beta < \alpha$}. 
$$
We can characterize $H_{\gamma}(0,T;\R)$ (see Appendix 5.1).
However, we do not need the characterization essentially in this 
article. The following lemma shows the fundamental 
property that $H_{\gamma}(0,T;X)$ is a Banach space. 
\begin{lemma}
\label{lem:Hgam}
$\Hgam$ is a Banach space with the norm $\Vert \cdot\Vert_{\Hgam}$.
\end{lemma}
\begin{proof}
We can readily prove that $\Hgam$ is a normed space.
We will prove the completeness: 
let $\{u_n\}_{n=1}^\infty$ be an arbitrary Cauchy sequence in 
$H_\gamma(0,T;X)$, which yields
$\lim_{n,m\to\infty}
\Vert u_n-u_m\Vert_{\Hgam} = 0$.  Then the definition of $J^{-\gamma}$ implies
$\lim_{n,m\to\infty} \Vert J^{-\gamma}u_n - J^{-\gamma}u_m\Vert
_{\LLLLLL} = 0$.  Since $\LLLLLL$ is complete, there exists a function  
$v_0 \in \LLLLLL$ such that 
$$
\lim_{n\to\infty} \Vert J^{-\gamma} u_n - v_0\Vert_{\LLLLLL}
= 0.
$$
The definition of the norm yields
$$
\Vert u_n - \JJJa v_0\Vert_{\Hgam}
= \Vert \JJJa (J^{-\gamma} u_n - v_0) \Vert_{\Hgam}
= \Vert J^{-\gamma} u_n - v_0\Vert_{\LLLLLL}.
$$
Therefore, $\lim_{n\to \infty} \Vert u_n - J^{\gamma}v_0\Vert_{\Hgam}
= 0$, and we obtain $J^{\gamma}v_0 \in J^{\gamma}L^2(0,T;X) = \Hgam$.
Thus, the proof of Lemma \ref{lem:Hgam} is complete.
\end{proof}

Now we define the fractional derivative by 
$$
\ppp_t^{\gamma}:= (J^{\gamma})^{-1}, \quad \DDD(\pppg) 
= H_{\gamma}(0,T;X).
$$

Throughout this article, we emphasize that $\frac{d}{dt}$ means the 
pointwise differentiation and $\ppp_t^1$ is an operator with the 
domain $H_1(0,T;X)$ where the operation is the same as 
$\frac{d}{dt}$.

Henceforth, we denote a norm equivalence by $\sim$. 
By the definition \eqref{intro:eq2-1} of $H_\gamma(0,T;X)$, we can 
readily see
\begin{theorem}
\label{thm:J2}
Let $\gamma > 0$. 

\noindent $\mathrm{(i)}$  
$J^{\gamma}: \LLLLLL \longrightarrow \Hgam$ 
is injective and surjective.

\noindent $\mathrm{(ii)}$ 
$\pppg J^{\gamma}u = u$ for all $u \in \LLLLLL$ and
$J^{\gamma}\pppg v = v$ for all $v \in \Hgam$.
\end{theorem}
Here we list properties of $\partial_t^{\alpha}$ as Theorems \ref{thm:J3} and \ref{thm:J4}.
\begin{theorem}
\label{thm:J3}
Let $\alpha>0$, $\beta\ge 0$. 

\noindent $\mathrm{(i)}$ $J^\alpha:$ $H_\beta(0,T;X)\longrightarrow H_{\alpha+\beta}(0,T;X)$ 
is injective and surjective, and \\
$\|J^\alpha u\|_{H_{\alpha+\beta}(0,T;X)} \sim \|u\|_{H_\beta(0,T;X)}$ for 
$u\in H_\beta(0,T;X)$. 

\noindent $\mathrm{(ii)}$ $\partial_t^\alpha:$ $H_{\alpha+\beta}(0,T;X)\longrightarrow 
H_{\beta}(0,T;X)$ is injective and surjective, and \\
$\|\partial_t^\alpha u\|_{H_{\beta}(0,T;X)} 
\sim \|u\|_{H_{\alpha+\beta}(0,T;X)}$ for $u\in H_{\alpha+\beta}(0,T;X)$. 
\end{theorem}

\begin{theorem}
\label{thm:J4}
$\mathrm{(i)}$ Let $\alpha,\beta>0$. 
Then $\partial_t^{\alpha+\beta} u = \partial_t^{\alpha}\partial_t^{\beta} u 
= \partial_t^{\beta}\partial_t^{\alpha} u$ for 
$u\in H_{\alpha+\beta}(0,T;X)$.

\noindent $\mathrm{(ii)}$ Let $0 < \beta \le \alpha$. Then 
$J^{\alpha-\beta} u = \partial_t^\beta J^\alpha u$ for 
$u\in L^2(0,T;X)$ and \\
$J^{\alpha-\beta} u = J^\alpha \partial_t^\beta u$ for 
$u\in H_\beta(0,T;X)$.

\noindent $\mathrm{(iii)}$ Let $0<\alpha<\beta$. Then 
$\partial_t^{\beta-\alpha} u = \partial_t^\beta J^\alpha u$ for
$u\in H_{\beta-\alpha}(0,T;X)$ and \\
$\partial_t^{\beta-\alpha} u  = J^\alpha \partial_t^\beta u$ for 
$u\in H_\beta(0,T;X)$.
\end{theorem}

These properties may be known in different contexts, but 
our treatments are based on the 
operator theory.  In particular, we attach domains whenever we consider 
$\pppa$. The proofs of Theorems \ref{thm:J3} and \ref{thm:J4} are straightforward by using
the definition \eqref{intro:eq2-1} and Lemma \ref{lem:J1}(iii), but we provide the proof of Theorem \ref{thm:J4}(ii) for example.
\begin{proof}
Let $u \in \LLLLLL$.  Then $J^{\beta}J^{\alpha-\beta}u
= J^{\beta+(\alpha-\beta)}u = J^{\alpha}u$ by Lemma \ref{lem:J1}(iii).
Since $J^{\alpha}u \in \Halp$, operating $J^{-\beta}$, we obtain
$J^{\alpha-\beta}u = J^{-\beta}J^{\alpha}u = \ppp_t^{\beta}J^{\alpha}u$,
which proves the first equation in (ii).  

Next, let $u \in H_{\beta}(0,T;X)$.
Then by \eqref{intro:eq2-1}, we can find some $w \in \LLLLLL$ such that 
$u = J^{\beta}w$. Hence
$J^{\alpha-\beta}u = J^{\alpha-\beta}J^{\beta}w = J^{(\alpha-\beta)+\beta}w 
= J^{\alpha}w$.  Since $u = J^{\beta}w$ means that $w = \ppp_t^{\beta}u$, 
we see that $J^{\alpha-\beta}u = J^{\alpha}\ppp_t^{\beta}u$.
Thus, the proof of Theorem \ref{thm:J4}(ii) is complete.
\end{proof}

Moreover, we have the following proposition, which indicates that the fractional derivative 
coincides with the Riemann-Liouville fractional derivative (see e.g., Podlubny \cite{P99}) 
for functions in $H_\gamma(0,T;X)$.   
\begin{proposition}
\label{prop:expression}
$\mathrm{(i)}$ $H_1(0,T;X) = \{ u \in H^1(0,T;X);\, u(0) = 0\}$. 

\noindent $\mathrm{(ii)}$ Let $\gamma = m-1 + \theta$ with $m\in \N$ and $0<\theta \le 1$. Then 
$$
\partial_t^\gamma u 
= \frac{d^m}{dt^m} J^{1-\theta} u, \quad u\in H_\gamma(0,T;X).
$$
In particular, 
\begin{equation}
\label{intro:eq2-2}
\partial_t^{m} u = J^{-m} u = \frac{d^m}{dt^m} u,\quad 
u\in H_m(0,T;X).
\end{equation}
\end{proposition}
\begin{proof}
(i) If $u\in H^1(0,T;X)$ and $u(0)=0$, then we can find a function 
$v\in L^2(0,T;X)$ such that $u(t) = \int_0^t v(s) ds = J^1 v(t)$ 
for all $0<t<T$. 
The existence of such $v$ can be verified, for example, by introducing $u_{\ep}
:= u*\varphi_{\ep}$ with the mollifier $\varphi_{\ep}$ (e.g., Adams \cite{AF75})
with $\ep>0$ and then letting $\ep \downarrow 0$.
Thus, $u\in H_1(0,T;X)$ by the definition \eqref{intro:eq2-1}. We immediately obtain the opposite inclusion by noting that the distributional derivative of $J^1 v$ is $v$ for any $v\in L^2(0,T;X)$. 

\noindent (ii) First, we verify \eqref{intro:eq2-2}. For any $u\in H_m(0,T;X)$, we find a function $v\in L^2(0,T;X)$ such that $u = J^m v$. 
By calculating the $m$-th distributional derivative of $J^m v$, we obtain $\frac{d^m}{dt^m} u = v$, which is $\partial_t^m u$ by the definition of the fractional derivative. 

Next, for any $u\in H_\gamma(0,T;X)$, there exists $w\in L^2(0,T;X)$ such that
$u=J^\gamma w$. Then $w = \ppp_t^{\gamma}u$.
Moreover, Lemma \ref{lem:J1}(iii) yields $J^{1-\theta}u = J^{1-\theta}J^{\gamma}w
= J^m w \in H_m(0,T;X)$, and by \eqref{intro:eq2-2} we see
$$
\frac{d^m}{dt^m} J^{1-\theta}u = J^{-m}J^{1-\theta}u
= J^{-m}J^m w = w.
$$
Thus, the proof of the proposition is complete.
\end{proof}

We conclude this subsection with 
a lemma which justifies the density argument for proving the coercivity
of the fractional derivatives (Appendix \ref{subsec:app2}), and  
is useful also for the interpretation of the initial condition 
in the initial-boundary value problem below.
\begin{lemma}
\label{lem:dense}
Let $X$ be a Banach space.

\noindent $\mathrm{(i)}$ Let $\alpha > \frac{1}{2}$. Then
$$
\Halp \subset C([0,T];X),
$$
and $u(0) = 0$ for $u\in \Halp$.

\noindent $\mathrm{(ii)}$ Let $\alpha > \frac{3}{2}$. Then
$$
\Halp \subset C^1([0,T];X),
$$
and $\frac{d}{dt} u(0) = 0$ for $u\in \Halp$.

\noindent $\mathrm{(iii)}$ For $0<\alpha<1$, the space
$$
\{ u\in C^1([0,T];X);\,  u(0)  = 0\},
$$
is dense in $\Halp$.
\end{lemma}
The proof can be found in Samko, Kilbas and Marichev \cite{SKM93} 
in the case of $X = \R$, and the proof in our case is essentially the same,
but for completeness we provide the proof in Appendix \ref{subsec:app2}.

\subsection{Statement of the main results}
\label{subsec:main-results}

Let $\Omega$ be a bounded domain in $\mathbb{R}^n$ with smooth boundary 
(e.g., of $\mathcal C^2$-class). 
Put $Q:= \Omega\times (0,T)$ and $\Sigma := \partial\Omega\times (0,T)$, 
where $T>0$ is arbitrarily fixed. 
In what follows, $A$ is defined by
$$
A (x,t)u := -\sum_{i,j=1}^n \partial_{i}\left(a_{ij}(x,t)\partial_{j}u\right) 
+ \sum_{j=1}^n b_j(x,t)\partial_j u + c(x,t) u,  \quad (x,t)\in \overline{Q},
$$
where $b_j, c\in W^{1,\infty}(0,T;L^\infty(\Omega))$, $1 \le j \le n$, 
$a_{ij}=a_{ji} \in W^{1,\infty}(Q)$, $1 \le i,j \le n$, satisfy
\begin{equation}
\label{assump:aij}
\sigma_0|\xi|^2 \le \sum_{i,j=1}^n a_{ij}(x,t)\xi_i\xi_j \le \sigma_1|\xi|^2, 
\quad (x,t)\in Q,\ \xi=(\xi_1,\ldots,\xi_n)\in \mathbb{R}^n,
\end{equation}
with some positive constants $\sigma_0,\sigma_1>0$. 

Identifying the dual space of $L^2(\OOO)$ with itself, we define a
Gel'fand triple: $H^1_0(\OOO) \subset L^2(\OOO) \subset 
(H^1_0(\OOO))'$: the dual space of $H^1_0(\OOO)$.  Henceforth we write
$H^{-1}(\OOO) = (H^1_0(\OOO))'$.

We formulate an initial-boundary value problem 
for the time-fractional diffusion-wave/wave equation ($1<\alpha\le 2$): 
\begin{equation}
\label{equ:pde2}
\left\{
\begin{aligned}
& \partial_t^\alpha (u - a_0 - t a_1) = -A(x,t) u(x,t) + F(x,t), &&\quad (x,t)\in \Omega\times (0,T),\\
& u(\cdot,t)\in H_0^1(\Omega), &&\quad t\in (0,T),\\
& u - a_0 - t a_1 \in H_\alpha(0,T;H^{-1}(\Omega)), 
\end{aligned}
\right.
\end{equation}
where $H_\alpha(0,T;H^{-1}(\OOO))$ 
and $\partial_t^\alpha$ denote the fractional Sobolev space
and the fractional derivative, respectively, which are introduced in the above 
subsections with $X=H^{-1}(\OOO)$. 
Here we impose the initial conditions 
by assuming $u - a_0 - t a_1 \in H_\alpha(0,T;H^{-1}(\Omega))$.  
In fact, in terms of Lemma \ref{lem:dense},
one can find that 
this is equivalent to $u(\cdot,0)=a_0$ and $\frac{\partial}{\partial t} u(\cdot,0)=a_1$ 
in $\Omega$ as long as the solution $u$ is sufficiently smooth. 

Now we are ready to state the main theorems, which are concerned 
with a weak solution. 
\begin{theorem}
\label{thm:pde1}
Let $F\in L^2(0,T;L^2(\Omega))$, $a_0\in H_0^1(\Omega)$ and $a_1\in L^2(\Omega)$. 
Then there exists a unique solution $u\in H^1(0,T;L^2(\Omega)) \cap L^\infty(0,T;H_0^1(\Omega))$ 
to \eqref{equ:pde2} satisfying $u - a_0 - t a_1 \in H_\alpha(0,T;H^{-1}(\Omega))$ and $\partial_t^{\alpha-1}(u - a_0)\in L^\infty(0,T;L^2(\Omega))$. 
Moreover, there exists a constant $C=C(\Omega,T,\alpha,\sigma_0,\sigma_1)>0$ such that
\begin{align}
\nonumber
&\quad \|u-a_0-t a_1\|_{H_\alpha(0,T;H^{-1}(\Omega))} 
+ \|\partial_t^{\alpha-1}(u-a_0)\|_{L^\infty(0,T;L^2(\Omega))}
+ \|u\|_{L^\infty(0,T;H_0^1(\Omega))} \\
\nonumber
&\quad + \|u\|_{H^1(0,T;L^2(\Omega))} \\
\label{esti:rev-pde1}
&\le C\left(\|a_0\|_{H_0^1(\Omega)} + \|a_1\|_{L^2(\Omega)} + \|F\|_{L^2(0,T;L^2(\Omega))}\right).
\end{align}
\end{theorem}
\begin{remark}
Here we should emphasize that Theorem \ref{thm:pde1} 
(as well as the following Theorem \ref{thm:pde2}) holds true 
for arbitrary $\alpha\in (1,2]$, that is, we can include 
the integer case of $\alpha=2$. Actually if we take $\alpha=2$, 
then Theorem \ref{thm:pde1} gives the well-known regularity 
for hyperbolic equations (e.g., Theorem 5 (p.410) 
in Section 7.2.3 of Evans \cite{E98}).  
Moreover, formally setting $\alpha=1$ in (5), we reach the classical 
regularity for parabolic equations, which coincides with   
e.g., Theorem 5 (p.382) in Section 7.1.3 of \cite{E98},  
except for the regularity of $L^2(0,T;H^2(\Omega))$. 
\end{remark}

If $a_0,a_1$ and $F$ are more regular, then we can gain the improved 
regularity as follows. 
\begin{theorem}
\label{thm:pde2}
Let $F-Aa_0\in H_1(0,T;L^2(\Omega))$, $a_0\in H^2(\Omega)\cap H_0^1(\Omega)$ and 
$a_1\in H_0^1(\Omega)$. 
Then the solution to \eqref{equ:pde2} satisfies
$u\in H^1(0,T;H_0^1(\Omega)) \cap L^\infty(0,T;H^2(\Omega)\cap H_0^1(\Omega))$,  
$\partial_t^\alpha (u - a_0 - t a_1) \in L^\infty(0,T;L^2(\Omega))$ and $\partial_t^{\alpha-1} (u - a_0) \in L^\infty(0,T;H_0^1(\Omega))$. 
Moreover, there exists a constant $C=C(\Omega,T,\alpha,\sigma_0,\sigma_1)>0$ such that
\begin{align*}
&\quad \|\partial_t^\alpha (u-a_0-t a_1)\|_{L^\infty(0,T;L^2(\Omega))} 
\!+\! \|\partial_t^{\alpha-1}(u-a_0)\|_{L^\infty(0,T;H_0^1(\Omega))}
\!+\! \|u\|_{L^\infty(0,T;H^2(\Omega))} \\
&\quad + \|u\|_{H^1(0,T;H_0^1(\Omega))} \\
&\le C\left(\|a_0\|_{H^2(\Omega)} + \|a_1\|_{H_0^1(\Omega)} + \|F\|_{H^1(0,T;L^2(\Omega))}\right).
\end{align*}
\end{theorem}
\begin{remark}
In Theorem \ref{thm:pde2}, we need the assumption $F-Aa_0\in H_1(0,T; L^2(\Omega))$ 
in order to obtain the improved regularity of the solution. 
In terms of \eqref{intro:eq2-2}, 
this assumption is equivalent to the following two assumptions: 
$$
\mathrm{(\romannumeral1)}\quad F\in H^1(0,T; L^2(\Omega)),\qquad\qquad 
\mathrm{(\romannumeral2)}\quad F(x,0) = A(x,0)a_0(x), \quad \mbox{for }\ x\in \Omega. 
$$
The former one (\romannumeral1) is a regularity assumption 
while the latter one (\romannumeral2) can be regarded as a compatible condition 
between the source term $F$ and the initial condition $a_0$. 
\end{remark}

The remaining part of this paper is organized as follows. 
In order to prove the main results, we first discuss 
the unique existence of the solution to a system of time-fractional 
ordinary differential equations in Section \ref{sec:ode}. 
After that, Sections \ref{sec:pde} and \ref{sec:pde2} are devoted to the proofs of 
Theorems \ref{thm:pde1} and \ref{thm:pde2}, respectively. 

\section{System of time-fractional ordinary differential equations}
\label{sec:ode}

For a single time-fractional ordinary differential equation with constant coefficients, 
the initial value problem was investigated by giving the explicit solution formula, 
see e.g., \cite{GKMR14}. 
In this section, for the sake of the proofs of Theorems \ref{thm:pde1} and \ref{thm:pde2},
we establish the well-posedness of the initial value problem for a system of 
time-fractional ordinary differential equations with the fractional order $\alpha\in (1,2]$. 

Let $N\in \N$. We discuss the solution $u=(u_1,\ldots,u_N)^T$ to the following 
initial value problem for a linear system of time-fractional ordinary differential equations:
\begin{equation}
\label{equ:ode2}
\left\{
\begin{aligned}
& \partial_t^\alpha \left(u - a_0 - t a_1\right) 
= P(t) u(t) + F(t),\quad 0<t<T,\\
& u - a_0 - t a_1 \in (H_\alpha(0,T))^N, 
\end{aligned}
\right.
\end{equation}
where $a_j = (a_{j,1}, \ldots, a_{j,N})^T$, $j=0,1$, $F=(f_1,\ldots,f_N)^T$ is given and 
$P=(p_{ij})_{i,j=1}^N$ is a given $N\times N$ matrix function. 
We have the following theorem. 
\begin{theorem}
\label{thm:ode2}
Let $\alpha\in (1,2]$, $a_0, a_1\in \mathbb{R}^N$, $F\in (L^2(0,T))^N$ and 
$P\in (L^\infty(0,T))^{N\times N}$. 
Then \eqref{equ:ode2} admits a unique solution 
$u\in (L^2(0,T))^N$ satisfying
\begin{equation}
\label{esti:ode2}
\left\|u - a_0 - t a_1\right\|_{(H_\alpha(0,T))^N} 
\le C\left(|a_0|_{\mathbb{R}^N} + |a_1|_{\mathbb{R}^N} + \|F\|_{(L^2(0,T))^N}\right),
\end{equation}
for some positive constant $C=C(\|P\|_{(L^\infty(0,T))^{N\times N}},N,T,\alpha)>0$. 
\end{theorem}
\begin{remark}
We can generalize Theorem \ref{thm:ode2} for all $\alpha>0$. 
Let $m-1<\alpha \le m$ with some $m\in \mathbb{N}$. We consider 
\begin{equation*}
\partial_t^\alpha \left(u - \sum_{j=0}^{m-1} \frac{t^j}{\Gamma(j+1)}a_j\right) 
= P(t) u(t) + F(t),\quad 0<t<T,
\end{equation*}
for $a_j\in \mathbb{R}^N$, $j=0,\ldots,m-1$. 
Then we can prove the unique existence of solution and a similar estimate as \eqref{esti:ode2} 
in the same way. 
\end{remark}
\begin{proof}
Let 
$$
v(t) = u(t) - a_0 - t a_1, \quad 
\widetilde F(t) = F(t) + P(t)a_0 + t P(t)a_1. 
$$
Then the proof is reduced to the verification that 
there exists a unique solution $v\in (H_\alpha(0,T))^N$ satisfying
$$
\partial_t^\alpha v(t) = P(t)v(t) + \widetilde F(t), \quad 0<t<T,
$$
for $\widetilde {F} \in (L^2(0,T))^N$.
By Theorem \ref{thm:J2}, we apply $J^\alpha$ on both sides of the above equation 
to obtain the integral equation:
\begin{equation}
\label{equ:integral}
v(t) = J^\alpha (Pv)(t) + J^\alpha \widetilde F(t),\quad 0<t<T.
\end{equation}
Henceforth by $H_{\gamma}(0,T)$ we denote $H_{\gamma}(0,T;\R)$.
 
Since $J^{\alpha}$ is an isomorphism from $L^2(0,T)$ to $H_\alpha(0,T)$ by the definition, 
it is sufficient to prove that there exists a unique solution to \eqref{equ:integral} in $(L^2(0,T))^N$. 
The unique solvability of the above integral equation can be 
obtained by the general theory on Volterra equations in some monographs 
(e.g., \cite{GLS90,P93}). 
However, in our case, we can take advantage of the compactness and provide 
a simple proof by the Fredholm alternative. 

Let $K: (L^2(0,T))^N \longrightarrow (L^2(0,T))^N$ be defined by
$$
K w := J^\alpha (Pw), \quad w\in (L^2(0,T))^N.
$$
According to the assumption that $P\in (L^\infty(0,T))^{N\times N}$, 
it is clear that $Pw\in (L^2(0,T))^N$. 
Thus, by Theorem \ref{thm:J2}, we have $J^\alpha (Pw)\in (H_\alpha(0,T))^N \hookrightarrow (L^2(0,T))^N$, 
which is compact.  
Then $K$ is a compact operator from $(L^2(0,T))^N$ to itself. 
By assuming $Kw = w$, we directly estimate 
\begin{align*}
|w(t)|_{\mathbb{R}^N} 
&= |Kw(t)|_{\mathbb{R}^N} = \left| \int_0^t \frac{(t-s)^{\alpha-1}}{\Gamma(\alpha)} P(s)w(s)ds \right|_{\mathbb{R}^N} \\
&\le \frac{\sqrt{N}}{\Gamma(\alpha)}\|P\|_{(L^\infty(0,T))^{N\times N}} 
\int_0^t (t-s)^{\alpha-1}|w(s)|_{\mathbb{R}^N}ds\\
&\le C \int_0^t (t-s)^{\alpha-1}|w(s)|_{\mathbb{R}^N}ds,\\
&\le CT^{\alpha-1} \int_0^t |w(s)|_{\mathbb{R}^N}ds,
\end{align*}
which together with Gr\"onwall's inequality ($\gamma=1$ in 
Lemma \ref{lem:Gronwall} in Section 5: Appendix) implies
\begin{equation}
\label{eq:w=0}
|w(t)|_{\mathbb{R}^N}=0,\quad 0<t<T.
\end{equation}
Therefore, by the Fredholm alternative, we reach that there exists a unique solution $v\in (L^2(0,T))^N$ 
to \eqref{equ:integral}. 
Moreover, we have
$$
v = J^\alpha (Pv + \widetilde F) \in (H_\alpha(0,T))^N.
$$
This proves the unique existence of the solution to \eqref{equ:ode2}. 

Next, we show the estimate \eqref{esti:ode2}. 
According to \eqref{equ:integral} and the estimation
$$
(t-s)^{\alpha-1} = (t-s) (t-s)^{\alpha-2} \le T (t-s)^{\alpha-2}, \quad 0 \le s < t \le T,
$$
we have
\begin{align*}
|v(t)|_{\mathbb{R}^N} 
&\le C \int_0^t (t-s)^{\alpha-1}|v(s)|_{\mathbb{R}^N}ds 
+ C \int_0^t (t-s)^{\alpha-1}|\widetilde F(s)|_{\mathbb{R}^N}ds\\
&\le C \int_0^t (t-s)^{\alpha-2}|v(s)|_{\mathbb{R}^N}ds 
+ C \int_0^t (t-s)^{\alpha-1}|\widetilde F(s)|_{\mathbb{R}^N}ds.
\end{align*}
By the generalized Gr\"onwall's inequality ($\gamma=\alpha-1$ in Lemma \ref{lem:Gronwall} in Section 5) and Fubini's theorem, we obtain
\begin{align*}
|v(t)|_{\mathbb{R}^N} 
&\le C \!\int_0^t \!(t-s)^{\alpha-1}|\widetilde F(s)|_{\mathbb{R}^N}ds 
\!+\! C \!\int_0^t \!(t-s)^{\alpha-2}\!\!\int_0^s \!(s-\tau)^{\alpha-1}|\widetilde F(\tau)|_{\mathbb{R}^N} 
d\tau ds\\
&\le C \int_0^t (t-s)^{\alpha-1}|\widetilde F(s)|_{\mathbb{R}^N}ds 
+ C\int_0^t (t-\tau)^{2\alpha-2}|\widetilde F(\tau)|_{\mathbb{R}^N} d\tau.
\end{align*}
Since $\alpha-1>0$, 
we take $L^2(0,T)$-norm of the above inequality and by Young's convolution inequality 
(e.g., \cite[Appendix A]{KRY20}), we obtain
$$
\|v\|_{(L^2(0,T))^N} \le C\|\widetilde F\|_{(L^2(0,T))^N}. 
$$
Hence together with the definition \eqref{intro:eq2-1} and \eqref{equ:integral}, 
this inequality yields
\begin{align*}
\|v\|_{(H_\alpha(0,T))^N} 
&= \|\partial_t^\alpha v\|_{(L^2(0,T))^N} \\
&\le C(\|Pv\|_{(L^2(0,T))^N} + \|\widetilde F\|_{(L^2(0,T))^N}) 
\le C\|\widetilde F\|_{(L^2(0,T))^N}.
\end{align*}
Therefore, we prove \eqref{esti:ode2} by combining the above inequality and 
the following estimate:
\begin{align*}
\|\widetilde F\|_{(L^2(0,T))^N} 
&\le \|F\|_{(L^2(0,T))^N} + \left\|Pa_0 + t Pa_1\right\|_{(L^2(0,T))^N}\\
&\le C\left(\|F\|_{(L^2(0,T))^N} + |a_0|_{\mathbb{R}^N} + |a_1|_{\mathbb{R}^N}\right).
\end{align*}
\end{proof}

We end up this section with the following result of an improved regularity. 
\begin{theorem}
\label{thm:ode3}
Let $\alpha\in (1,2]$, $a_0, a_1\in \mathbb{R}^N$ and $m-1<\beta\le m$ for some
$m\in \N$.  
We assume $P\in (W^{m,\infty}(0,T))^{N\times N}$ 
and 
$$
F+ Pa_0 + t Pa_1
\in (H_\beta(0,T))^N. 
$$
Then the solution $u$ to \eqref{equ:ode2} satisfies
$$
u - a_0 - t a_1
\in (H_{\alpha+\beta}(0,T))^N.
$$
\end{theorem}

\begin{proof}
By the proof of Theorem \ref{thm:ode2}, the solution to \eqref{equ:ode2} admits 
\begin{equation*}
u - a_0 - t a_1 = J^\alpha \left(Pu + F \right).
\end{equation*}
By setting 
$$
v(t) = u(t) - a_0 - t a_1,
\quad t\in [0,T], 
$$
we rewrite the above equation by 
\begin{equation}
\label{eq:ode-im}
v = J^\alpha \left(Pv\right) 
+ J^\alpha \left(F + Pa_0 + t Pa_1 \right).
\end{equation}
Since $v\in (H_\alpha(0,T))^N$ and $P\in (W^{m,\infty}(0,T))^{N\times N}$,
by $\beta\le m< \beta+1$, 
we have $Pv\in (H_\beta(0,T))^N$ if $\beta\le \alpha$ 
and $Pv\in (H_\alpha(0,T))^N$ if $\beta> \alpha$, 
that is, $Pv\in (H_{\min\{\alpha,\beta\}}(0,T))^N$. 
Recalling the assumption $F+ Pa_0 + t Pa_1\in (H_\beta(0,T))^N$,  
the equation \eqref{eq:ode-im} yields 
$v\in (H_{\min\{2\alpha,\alpha+\beta\}}(0,T))^N$. 
If $2\alpha\ge \alpha+\beta$, then it is done. Otherwise, we have $Pv\in (H_{\min\{2\alpha,\beta\}}(0,T))^N$. 
Then again by \eqref{eq:ode-im}, we further obtain $v\in (H_{\min\{3\alpha,\alpha+\beta\}}(0,T))^N$. 
We end up with $v\in (H_{\alpha+\beta}(0,T))^N$ by iterating for $k_0$ times until $k_0\alpha\ge \alpha+\beta$. 
\end{proof}

We finish this section with a special case $\beta=1$ of Theorem \ref{thm:ode3},
which will be used in Section \ref{sec:pde2}. 
\begin{corollary}
\label{coro:ode4}
Let $a_0, a_1\in \mathbb{R}^N$ 
and $P\in (W^{1,\infty}(0,T))^{N\times N}$. Moreover, we assume 
$$
F+ Pa_0 \in (H_1(0,T))^N. 
$$
Then the solution $u$ to \eqref{equ:ode2} satisfies
$$
u - a_0 - t a_1
\in (H_{\alpha+1}(0,T))^N.
$$
\end{corollary}
This result follows immediately from Theorem \ref{thm:ode3} by noting the fact that 
$F+ Pa_0 \in (H_1(0,T))^N$ is equivalent to $F+ Pa_0 + t Pa_1\in (H_1(0,T))^N$ provided that $P\in (W^{1,\infty}(0,T))^{N\times N}$. 
\section{Proof of Theorem \ref{thm:pde1}}
\label{sec:pde}

In this and the next sections, we apply the Galerkin approximation to 
prove Theorems \ref{thm:pde1} and \ref{thm:pde2} respectively. 
We divide the proof into four steps. 

\noindent {\bf Step 1. Construction of the approximate solutions} 

We introduce $0<\lambda_1\le \lambda_2 \le \ldots$ as the eigenvalues of $-\Delta$ with 
the homogeneous Dirichlet boundary condition, where we number $\lambda_k$ 
according to the multiplicity, that is, if the multiplicity of $\lambda_k$ is
$m$, then $\lambda_k$ appears repeatedly $m$-times in this sequence. 
By $\varphi_k\in H^2(\Omega)\cap H_0^1(\Omega)$, $k\in \N$, we denote 
the corresponding eigenfunction of $-\Delta$ for $\lambda_k$. 
We can choose $\varphi_k$ such that $\{\varphi_k\}_{k=1}^\infty$ forms 
an orthonormal basis in $L^2(\Omega)$, that is, we have
\begin{equation}
(\varphi_k,\varphi_\ell)_{L^2(\Omega)} = \delta_{k\ell}, \quad 
_{H^{-1}(\Omega)}\langle \varphi_k,\varphi_\ell \rangle_{H_0^1(\Omega)} = \delta_{k\ell},
\end{equation}
where we set
$$
\delta_{k\ell} := 
\left\{
\begin{aligned}
& 1,\quad && k=\ell,\\
& 0,\quad && k\not=\ell.
\end{aligned}
\right.
$$
We note by $\varphi_k, \varphi_{\ell} \in H^1_0(\OOO)$ that 
$(\varphi_k,\varphi_\ell)_{L^2(\Omega)} = \, 
_{H^{-1}(\Omega)}\langle \varphi_k,\varphi_\ell \rangle_{H_0^1(\Omega)}$.

For arbitrarily fixed $N\in \N$, we seek for the approximate solutions
in the form 
$$
u_N(x,t) = \sum_{k=1}^N p_k^N(t)\varphi_k(x),
$$
satisfying
\begin{equation}
\label{equ:pde2b}
\left\{
\begin{aligned}
& \partial_t^\alpha (u_N - a_{0,N} - t a_{1,N}) 
= -A(x,t) u_N(x,t) + F_N(x,t), &&\quad (x,t)\in \Omega\times (0,T),\\
& u_N - a_{0,N} - t a_{1,N} \in H_\alpha(0,T;L^2(\Omega)), 
\end{aligned}
\right.
\end{equation}
where
\begin{align*}
&a_{j,N} = \sum_{k=1}^N a_k^j \varphi_k, \quad 
a_k^j := (a_j,\varphi_k)_{L^2(\Omega)}, \quad j=0,1,\\
&F_N(t) = \sum_{k=1}^N f_k(t)\varphi_k, \quad 
f_k(t) := _{H^{-1}(\Omega)}\langle F(t),\varphi_k \rangle_{H_0^1(\Omega)}, \quad 0<t<T.
\end{align*}
For each $\ell = 1,\ldots,N$, we multiply $\varphi_\ell$ on both sides of \eqref{equ:pde2b} and 
take integral over $\Omega$. Then we obtain
$$
\left\{
\begin{aligned}
& \partial_t^\alpha (p_\ell^N - a_\ell^0 - t a_\ell^1) 
= \sum_{k=1}^N p_k^N \left(-A(t)\varphi_k,\varphi_\ell\right)_{L^2(\Omega)} + f_\ell(t), 
&&\quad t\in (0,T),\\
& p_\ell^N - a_\ell^0 - t a_\ell^1 \in H_\alpha(0,T),  &&
\quad \ell = 1,\ldots,N.
\end{aligned}
\right.
$$
We define $Q(t)=(q_{\ell k}(t))_{k,\ell=1}^N$ by 
\begin{equation}
\label{def:Q}
q_{\ell k}(t) := \left(-A(t)\varphi_k,\varphi_\ell\right)_{L^2(\Omega)}, \quad 0<t<T.
\end{equation}
By the regularity assumptions on the coefficients of $A$ and $F\in L^2(0,T;L^2(\Omega))$, 
we have 
$$
Q\in (L^\infty(0,T))^{N\times N},\quad \mbox{and} \quad
f = (f_1,\ldots,f_N)^T\in (L^2(0,T))^N. 
$$
Then we rewrite the above equations by setting $a^j = (a_1^j,\ldots,a_N^j)^T$ for $j=0,1$, 
$p^N = (p_1^N,\ldots,p_N^N)^T$: 
\begin{equation}
\label{equ:pde2c}
\left\{
\begin{aligned}
& \partial_t^\alpha (p^N - a^0 - t a^1) 
= Q(t)p^N(t) + f(t), 
&&\quad t\in (0,T),\\
& p^N - a^0 - t a^1 \in (H_\alpha(0,T))^N.
\end{aligned}
\right.
\end{equation}
Therefore, by Theorem \ref{thm:ode2}, we find a unique solution $p^N\in (L^2(0,T))^N$ 
satisfying $p^N - a^0 - t a^1\in (H_\alpha(0,T))^N$, and we construct the approximate solution 
$u_N$ which is the solution to \eqref{equ:pde2b} satisfying
$$
u_N\in L^2(0,T;H^2\cap H_0^1(\Omega)), \quad 
u_N - a_{0,N} - t a_{1,N}\in H_\alpha(0,T;H^2\cap H_0^1(\Omega)). 
$$
By varying $N$ in $\N$, we obtain the approximate solutions $\{u_N\}_{N=1}^\infty$. 

\noindent {\bf Step 2. A priori estimate of the approximate solutions}

Henceforth, we may write 
$u_N(\cdot,t)=u_N(t)=u_N$, $a_{ij}(x,t) = a_{ij}(t) = a_{ij}$,
etc. according to the convenience if there is no fear of confusion. 

We multiply $\frac{\ppp u_N}{\ppp t}$ on both sides of the equation \eqref{equ:pde2b} and then 
integrate over $\Omega$. 
Moreover, we apply the fractional integral operator $J^{\alpha-1}$ and obtain 
\begin{align}
\nonumber
&\quad J^{\alpha-1}\int_\Omega \ppp_t^\alpha (u_N-a_{0,N}-ta_{1,N}) 
\pppppp (u_N-a_{0,N}-ta_{1,N}) dx \\
\nonumber
&\quad + J^{\alpha-1}\int_\Omega \partial_t^\alpha (u_N-a_{0,N}-ta_{1,N}) 
a_{1,N} dx\\
\label{eq:ap1}
&= J^{\alpha-1}\int_\Omega -Au_N \pun dx 
+ J^{\alpha-1} \int_\Omega F_N \pun dx. 
\end{align}
As we find $u_N-a_{0,N}-t a_{1,N}\in H_\alpha(0,T;H^2(\Omega)\cap H_0^1(\Omega))$ 
$\subset H_1(0,T;H^2(\Omega)\cap H_0^1(\Omega))$, 
we have by Theorem \ref{thm:J4}(i) and Proposition \ref{prop:expression} that 
\begin{align*}
&\frac{\partial}{\partial t} (u_N-a_{0,N}-ta_{1,N}) = \partial_t^1 (u_N-a_{0,N}-ta_{1,N}), \quad \mbox{and}\quad\\
&\partial_t^\alpha (u_N\!-\!a_{0,N}\!-\!ta_{1,N}) 
= \partial_t^{\alpha-1}\partial_t^1 (u_N\!-\!a_{0,N}\!-\!ta_{1,N}) 
= \partial_t^{\alpha-1}\frac{\partial}{\partial t}(u_N\!-\!a_{0,N}\!-\!ta_{1,N}). 
\end{align*}
Then, we set  
$$
v(t) := \partial_t^1 (u_N-a_{0,N}-ta_{1,N}) = \frac{\partial}{\partial t}(u_N-a_{0,N}-ta_{1,N})
= \frac{\partial}{\partial t} u_N - a_{1,N},
$$
and see that $v\in H_{\alpha-1}(0,T;H^2(\Omega)\cap H_0^1(\Omega))$ and
$u_N-a_{0,N} - ta_{1,N} = J^1v$ by Theorem \ref{thm:J3}.
Hence, by Theorem \ref{thm:J4}(ii), we employ the coercivity inequality (Lemma \ref{lem:coer} in Appendix) 
to estimate the first term 
on the left-hand side of \eqref{eq:ap1}: 
\begin{align*}
&\quad J^{\alpha-1}\left(
\int_\Omega \partial_t^\alpha (u_N-a_{0,N}-ta_{1,N})(x,t) 
\pppppp (u_N-a_{0,N}-ta_{1,N})(x,t) dx\right) \\
&= (J^{\alpha-1}(\pppa J^1 v(\cdot,t), v(\cdot,t))_{L^2(\OOO)})(t)
= (J^{\alpha-1} (\partial_t^{\alpha-1} v(\cdot,t), v(\cdot,t))_{L^2(\Omega)})
(t)                             \\
&\ge \frac12 \|v(\cdot,t)\|_{L^2(\Omega)}^2 
= \frac12 \left\| \pppppp u_N(\cdot,t)-a_{1,N}\right\|_{L^2(\Omega)}^2\\
&\ge \frac14 \left\|\pppppp u_N(\cdot,t)\right\|_{L^2(\Omega)}^2 
- \frac12 \|a_{1,N}\|_{L^2(\Omega)}^2.
\end{align*}
Here in the last line we used the estimate $(A-B)^2\ge \frac12 A^2 - B^2$. 
For the second term on the left-hand side of \eqref{eq:ap1}, 
by Fubini's theorem and Theorem \ref{thm:J4}(i), we have
\begin{align*}
J^{\alpha-1}\int_\Omega \partial_t^\alpha (u_N-a_{0,N}-ta_{1,N}) a_{1,N} dx
&= \int_\Omega (J^{\alpha-1} \partial_t^{\alpha-1} v(x,t)) a_{1,N} dx\\
&= \left(\frac{\partial}{\partial t} u_N(\cdot,t),a_{1,N}\right)_{L^2(\Omega)} 
- \|a_{1,N}\|_{L^2(\Omega)}^2\\
&\ge -\frac18 \left\|\frac{\partial}{\partial t} u_N(\cdot,t)\right\|_{L^2(\Omega)}^2 
- 3\|a_{1,N}\|_{L^2(\Omega)}^2.
\end{align*}
Here in the last line we used the estimate $AB=(\frac12 A) (2B)\ge -\frac18 A^2 - 2 B^2$. 

On the other hand, we estimate the first term on the right-hand side of \eqref{eq:ap1} 
by integration by parts, and we have
\begin{align*}
&\quad J^{\alpha-1}\int_\Omega -Au_N \pppppp u_N dx \\
&= - J^{\alpha-1}\int_\Omega \sum_{i,j=1}^n a_{ij} (\partial_i u_N)
   \pppppp (\partial_j u_N) dx 
+ J^{\alpha-1}\int_\Omega \Big(\sum_{j=1}^n b_j \partial_j u_N 
+ c u_N\Big)\frac{\ppp u_N}{\ppp t} dx\\
&= -\frac12 J^{\alpha-1} \dddddd \int_\Omega \sum_{i,j=1}^n a_{ij} 
(\partial_i u_N) \partial_j u_N dx\\ 
&\quad + J^{\alpha-1}\int_\Omega \Big(\frac12 \sum_{i,j=1}^n 
\frac{\ppp a_{ij}}{\ppp t} 
(\partial_i u_N) \partial_j u_N
+ \sum_{j=1}^n b_j (\partial_j u_N) \frac{\ppp u_N}{\ppp t} 
+ c u_N \frac{\ppp u_N}{\ppp t}\Big) dx\\
&\le -\frac12 J^{\alpha-1} \dddddd\int_\Omega \sum_{i,j=1}^n a_{ij} 
(\partial_i u_N) \partial_j u_N dx 
\!+\! CJ^{\alpha-1}\|\nabla u_N(t)\|_{L^2(\Omega)}^2 \\
&\quad + \frac12 J^{\alpha-1}\left\|\frac{\partial u_N}{\partial t} (\cdot,t) \right\|_{L^2(\Omega)}^2.
\end{align*}
Moreover, we have
$$
J^{\alpha-1} \int_\Omega F_N \frac{\ppp u_N}{\ppp t} dx 
\le \frac12 J^{\alpha-1}\|F_N(\cdot,t)\|_{L^2(\Omega)}^2 
+ \frac12 J^{\alpha-1}\left\Vert \frac{\partial u_N}{\partial t} (\cdot,t)\right\Vert
_{L^2(\Omega)}^2.
$$
Therefore, combining the above four estimates and \eqref{eq:ap1} yields
\begin{align*}
&\quad \frac18 \left\|\frac{\partial u_N}{\partial t} (\cdot,t)\right\|_{L^2(\Omega)}^2 
+ \frac12 J^{\alpha-1} \dddddd \int_\Omega \sum_{i,j=1}^n a_{ij} 
(\partial_i u_N) \partial_j u_N dx\\
&\le \frac72\|a_{1,N}\|_{L^2(\Omega)}^2 
+ \frac12 J^{\alpha-1}\|F_N(\cdot,t)\|_{L^2(\Omega)}^2 
+ J^{\alpha-1}\left\|\frac{\partial u_N}{\partial t} (\cdot,t)\right\|_{L^2(\Omega)}^2 \\
&\quad + CJ^{\alpha-1}\|\nabla u_N(\cdot,t)\|_{L^2(\Omega)}^2,
\end{align*}
for $0<t<T$. 
Next, we apply the fractional integral operator $J^{2-\alpha}$ 
(it should be understood as the identity operator if $\alpha=2$) on both sides above. 
By noting 
$$
J^{2-\alpha}J^{\alpha-1} v = J^{\alpha-1}J^{2-\alpha} v = J^1 v = \int_0^t v(s)ds \quad 
\mbox{for } v\in L^2(0,T;\mathbb{R}),
$$ 
we obtain 
\begin{align}
\nonumber
&\quad \frac18 J^{2-\alpha} \left\| \frac{\partial u_N}{\partial t} (\cdot,t)
\right\|_{L^2(\Omega)}^2 
+ \frac12 \int_\Omega \sum_{i,j=1}^n a_{ij}(x,t) \partial_i u_N(x,t) 
\partial_j u_N(x,t) dx\\
\nonumber
&\quad - \frac12 \int_\Omega \sum_{i,j=1}^n a_{ij}(x,0) 
\partial_i u_N(x,0) \partial_j u_N(x,0) dx\\
\nonumber
&\le \frac72\frac{t^{2-\alpha}}{\Gamma(3-\alpha)}\|a_{1,N}\|_{L^2(\Omega)}^2 
\!+\! \frac12 J^1\|F_N(\cdot,t)\|_{L^2(\Omega)}^2 
\!+\! J^1\left\| \frac{\partial u_N}{\partial t} (\cdot,t)\right\|_{L^2(\Omega)}^2 \\
\label{eq:ap2}
&\quad + \! CJ^1\|\nabla u_N(\cdot,t)\|_{L^2(\Omega)}^2,
\end{align}
for $0<t<T$. By the assumption \eqref{assump:aij} on $a_{ij}$, we have
\begin{align*}
&\frac12 \int_\Omega \sum_{i,j=1}^n a_{ij}(x,t) 
\partial_i u_N(x,t) \partial_j u_N(x,t) dx 
\ge \frac{\sigma_0}{2}\|\nabla u_N(\cdot,t)\|_{L^2(\Omega)}^2,\quad \mbox{and } \\
&\frac12 \int_\Omega \sum_{i,j=1}^n a_{ij}(x,0) \partial_i u_N(x,0) 
\partial_j u_N(x,0) dx
\le \frac{\sigma_1}{2}\|\nabla u_N(\cdot,0)\|_{L^2(\Omega)}^2.
\end{align*}
Moreover, by Lemma \ref{lem:dense}, $u_N-a_{0,N}-t a_{1,N}\in H_\alpha(0,T;H^2(\Omega)\cap H_0^1(\Omega))$ implies $u_N(\cdot,0) = a_{0,N}$, and hence
$$
\|\nabla u_N(\cdot,0)\|_{L^2(\Omega)}^2 \sim \|u_N(\cdot,0)\|_{H_0^1(\Omega)}^2 
= \|a_{0,N}\|_{H_0^1(\Omega)}^2. 
$$
Thus, we put the above estimates into \eqref{eq:ap2}, and we obtain
\begin{align}
\nonumber
&\quad J^{2-\alpha} \left\| \frac{\partial u_N}{\partial t} (\cdot,t)\right\|_{L^2(\Omega)}^2 
+ \|\nabla u_N(\cdot,t)\|_{L^2(\Omega)}^2\\
\nonumber
&\le C \left(\|a_{1,N}\|_{L^2(\Omega)}^2 + \|a_{0,N}\|_{H_0^1(\Omega)}^2 
+ J^1\|F_N(\cdot,t)\|_{L^2(\Omega)}^2\right) \\
\label{eq:ap3}
&\quad + CJ^1\left\| \frac{\partial u_N}{\partial t} (\cdot,t)\right\|_{L^2(\Omega)}^2 
+ CJ^1\|\nabla u_N(\cdot,t)\|_{L^2(\Omega)}^2,
\end{align}
for $0<t<T$. 
Furthermore, by the property of the Riemann-Liouville fractional integral operator, we note that 
$J^1 \left(\left\| \frac{du_N}{dt}\right\|_{L^2(\Omega)}^2\right)(t) 
= J^{\alpha-1}J^{2-\alpha}
\left(\left\| \frac{du_N}{dt}\right\|_{L^2(\Omega)}^2\right)(t)$ 
and
\begin{align*}
J^1(\|\nabla u_N\|_{L^2(\Omega)}^2)(t) 
&= \int_0^t \|\nabla u_N(\cdot,s)\|_{L^2(\Omega)}^2 ds \\
&= \int_0^t \Gamma(\alpha-1)(t-s)^{2-\alpha} \frac{(t-s)^{\alpha-2}}{\Gamma(\alpha-1)}\|\nabla u_N(\cdot,s)\|_{L^2(\Omega)}^2 ds\\
&\le \Gamma(\alpha-1)T^{2-\alpha} J^{\alpha-1}\|\nabla u_N(\cdot,t)\|
_{L^2(\Omega)}^2, \quad 0<t<T. 
\end{align*}
Then by setting 
$$
w(t) := \left( J^{2-\alpha}\left\| \frac{du_N}{dt} \right\|_{L^2(\Omega)}^2
\right)(t) + \|\nabla u_N(\cdot,t)\|_{L^2(\Omega)}^2 \ge 0,
$$
we rewrite \eqref{eq:ap3} by 
$$
w(t) \le C\left(\|a_{0,N}\|_{H_0^1(\Omega)}^2 + \|a_{1,N}\|_{L^2(\Omega)}^2
+ \|F_N\|_{L^2(0,t;L^2(\Omega))}^2\right) + C J^{\alpha-1} w(t),
$$
for $0<t<T$. 
Therefore, by the generalized Gr\"onwall's inequality ($\gamma=\alpha-1$ in Lemma \ref{lem:Gronwall} in Section 5) , we obtain
$$
w(t) \le C\left(\|a_{0,N}\|_{H_0^1(\Omega)}^2 + \|a_{1,N}\|_{L^2(\Omega)}^2
+ \|F_N\|_{L^2(0,t;L^2(\Omega))}^2\right),\quad 0<t<T,
$$
with a new constant $C>0$. 
By the definition of $a_{j,N}$, $j=0,1$ and $F_N$, we see that
\begin{align*}
&\|a_{0,N}\|_{H_0^1(\Omega)}^2 \le \|a_0\|_{H_0^1(\Omega)}^2, \quad
\|a_{1,N}\|_{L^2(\Omega)}^2 \le \|a_1\|_{L^2(\Omega)}^2,\\
&\|F_N(\cdot,t)\|_{L^2(\Omega)}^2 \le \|F(\cdot,t)\|_{L^2(\Omega)}^2,\ 0<t<T,
\end{align*}
and hence we have
\begin{align}
\nonumber
&\quad J^{2-\alpha}\left(\left\| \frac{du_N}{dt}\right\|_{L^2(\Omega)}^2
\right)(t) 
+ \|\nabla u_N(\cdot,t)\|_{L^2(\Omega)}^2 \\
\label{eq:ap4-}
&\le C\left(\|a_0\|_{H_0^1(\Omega)}^2 + \|a_1\|_{L^2(\Omega)}^2 + \|F\|_{L^2(0,t;L^2(\Omega))}^2\right),
\end{align}
for $0<t<T$. 
Moreover, we have the following estimate by using H\"older's inequality and \eqref{eq:ap4-} if $1<\alpha<2$: 
\begin{align}
\nonumber
&\quad \left\|\left(J^{2-\alpha} \frac{\partial u_N}{\partial t}\right)(\cdot,t)\right\|_{L^2(\Omega)}^2 \\
\nonumber
&= \int_\Omega \left(\int_0^t \frac{(t-s)^{1-\alpha}}{\Gamma(2-\alpha)}
\frac{\ppp u_N}{\ppp s}(x,s)ds\right)^2 dx\\
\nonumber
&\le \int_\Omega \left(\int_0^t \frac{(t-s)^{1-\alpha}}{\Gamma(2-\alpha)}ds
\right) 
\left(\int_0^t \frac{(t-s)^{1-\alpha}}{\Gamma(2-\alpha)}
\left| \frac{\ppp u_N}{\ppp s}(x,s)\right|^2 ds\right) dx\\
\nonumber
&= \frac{t^{2-\alpha}}{\Gamma(3-\alpha)} J^{2-\alpha}
 \left(\left\| \frac{du_N}{dt} \right\|_{L^2(\Omega)}^2\right)(t)\\
\label{eq:ap4-2}
&\le C\left(\|a_0\|_{H_0^1(\Omega)}^2 + \|a_1\|_{L^2(\Omega)}^2 + \|F\|_{L^2(0,t;L^2(\Omega))}^2\right), 
\end{align}
for $0<t<T$. 
For $\alpha=2$, the above estimate is trivial by using \eqref{eq:ap4-}, that is, 
\begin{align*}
\left\|\left(J^{2-\alpha} \frac{\partial u_N}{\partial t}\right)(\cdot,t)\right\|_{L^2(\Omega)}^2 
&= \left(\left\| \frac{du_N}{dt} \right\|_{L^2(\Omega)}^2\right)(t) \\
&\le C\left(\|a_0\|_{H_0^1(\Omega)}^2 + \|a_1\|_{L^2(\Omega)}^2 + \|F\|_{L^2(0,t;L^2(\Omega))}^2\right).
\end{align*}
Since $u_N-a_{0,N}-t a_{1,N}\in H_\alpha(0,T;L^2(\Omega)) \subset H_{\alpha-1}(0,T;L^2(\Omega))$ and $t a_{1,N}\in H_{\alpha-1}(0,T;L^2(\Omega))$ imply $u_N-a_{0,N}\in H_{\alpha-1}(0,T;L^2(\Omega))$, by Theorem \ref{thm:J4}(iii), we have the equality
\begin{align}
\nonumber
\|\partial_t^{\alpha-1} (u_N-a_{0,N})(\cdot,t)\|_{L^2(\Omega)}^2 
&= \left\|J^{2-\alpha}\frac{\partial}{\partial t} (u_N-a_{0,N})(\cdot,t)\right\|_{L^2(\Omega)}^2 \\
\label{eq:ap4-3}
&= \left\|J^{2-\alpha}\frac{\partial u_N}{\partial t} (\cdot,t)\right\|_{L^2(\Omega)}^2.
\end{align}
Then we insert \eqref{eq:ap4-2} into the right-hand side of \eqref{eq:ap4-3} and use \eqref{eq:ap4-} again to obtain
\begin{align}
\nonumber
&\quad J^{2-\alpha}\left(
\left\| \frac{d u_N}{d t}\right\|_{L^2(\Omega)}^2\right)(t) 
+ \|\nabla u_N(\cdot,t)\|_{L^2(\Omega)}^2 
+ \|\partial_t^{\alpha-1} (u_N-a_{0,N})(\cdot,t)\|_{L^2(\Omega)}^2\\
\label{eq:ap4}
&\le C\left(\|a_0\|_{H_0^1(\Omega)}^2 + \|a_1\|_{L^2(\Omega)}^2 + \|F\|_{L^2(0,t;L^2(\Omega))}^2\right),
\end{align}
for $0<t<T$. 
By applying $J^{\alpha-1}$ on both sides of \eqref{eq:ap4} and then taking $t=T$, we use Young's convolution inequality 
(e.g., \cite[Appendix A]{KRY20}) for the third term on the right-hand side and conclude
\begin{equation}
\label{eq:ap5}
\left\| \frac{\partial u_N}{\partial t}\right\|_{L^2(0,T;L^2(\Omega))}^2 
\le C\left(\|a_0\|_{H_0^1(\Omega)}^2 + \|a_1\|_{L^2(\Omega)}^2 + \|F\|_{L^2(0,T;L^2(\Omega))}^2\right).
\end{equation}

Next, we estimate $\|u_N-a_{0,N}-t a_{1,N}\|_{H_\alpha(0,T;H^{-1}(\Omega))}$. 
For any $\psi\in H_0^1(\Omega)$, by the equation \eqref{equ:pde2b}, we have
\begin{align*}
&\quad \left|{_{H^{-1}(\Omega)}\langle \partial_t^\alpha (u_N-a_{0,N}-ta_{1,N}),\psi 
\rangle_{H_0^1(\Omega)}}\right| \\
&\le \left|{_{H^{-1}(\Omega)}\langle -Au_N(\cdot,t),\psi \rangle_{H_0^1(\Omega)}}\right|
+ \left|{_{H^{-1}(\Omega)}\langle F_N(\cdot,t),\psi \rangle_{H_0^1(\Omega)}}\right|\\
&\le C\left(\|\nabla u_N(\cdot,t)\|_{L^2(\Omega)} + \|F_N(\cdot,t)\|_{H^{-1}(\Omega)}\right) 
\|\psi\|_{H_0^1(\Omega)}^2.
\end{align*}
Therefore, by $\|F_N(\cdot,t)\|_{H^{-1}(\Omega)} \le C\|F_N(\cdot,t)\|_{L^2(\Omega)}$ 
$\le C\|F(\cdot,t)\|_{L^2(\Omega)}$ for $0<t<T$ and the above inequality, 
we find that $u_N-a_{0,N}-t a_{1,N}\in H_\alpha(0,T;H^{-1}(\Omega))$ with
\begin{align}
\nonumber
\|u_N\!-\!a_{0,N}\!-\!t a_{1,N}\|_{H_\alpha(0,T;H^{-1}(\Omega))} 
&\le C\|\partial_t^\alpha (u_N\!-\!a_{0,N}\!-\!t a_{1,N})\|_{L^2(0,T;H^{-1}(\Omega))} \\
\label{eq:ap6}
&\le C\left(\|\nabla u_N\|_{L^2(0,T;L^2(\Omega))} \!+\! \|F\|_{L^2(0,T;L^2(\Omega))}\right).
\end{align}
To sum up, we combine the estimates \eqref{eq:ap4}--\eqref{eq:ap6} and obtain the a priori estimate: 
\begin{align}
\nonumber
&\quad \|u_N\|_{L^\infty(0,T;H_0^1(\Omega))} + \|\partial_t^{\alpha-1}(u_N-a_{0,N})\|_{L^\infty(0,T;L^2(\Omega))} 
+ \left\| \frac{\partial u_N}{\partial t}\right\|_{L^2(0,T;L^2(\Omega))} \\
\nonumber
&\quad + \|u_N-a_{0,N}-t a_{1,N}\|_{H_\alpha(0,T;H^{-1}(\Omega))} \\
\label{eq:ap7}
&\le C\left(\|a_0\|_{H_0^1(\Omega)} + \|a_1\|_{L^2(\Omega)} + \|F\|_{L^2(0,T;L^2(\Omega))}\right).
\end{align}

\noindent {\bf Step 3. Existence of solution by taking a limit of a weakly convergent subsequence}

In the previous step, since the constant $C>0$ is independent of $N$, 
we derived a uniform estimate of some norms of the approximate solutions. 
By \eqref{eq:ap7}, in particular, we see that
\begin{align*}
&\quad \|u_N\|_{L^2(0,T;H_0^1(\Omega))} + \|u_N\|_{H^1(0,T;L^2(\Omega))} 
+ \|u_N-a_{0,N}-t a_{1,N}\|_{H_\alpha(0,T;H^{-1}(\Omega))} \\
&\le C\left(\|a_0\|_{H_0^1(\Omega)} + \|a_1\|_{L^2(\Omega)} + \|F\|_{L^2(0,T;L^2(\Omega))}\right).
\end{align*}
Therefore, the sequence $\{u_N\}_{N=1}^\infty$ is bounded in 
$L^2(0,T;H_0^1(\Omega)) \cap H^1(0,T;L^2(\Omega))$, and
$\{u_N-a_{0,N}-t a_{1,N}\}_{N=1}^\infty$ is bounded in $H_\alpha(0,T;H^{-1}(\Omega))$. 
Then there exist functions  
$u\in L^2(0,T;H_0^1(\Omega))\cap H^1(0,T;L^2(\Omega))$, 
$\widetilde u\in H_\alpha(0,T;H^{-1}(\Omega))$
and a subsequence $\{u_{N^\prime}\}$ of $\{u_N\}_{N=1}^\infty$ such that
\begin{equation}
\label{eq:convergence}
\left\{
\begin{aligned}
&u_{N^\prime} \rightharpoonup u 
&&\quad \mbox{weakly in}\ L^2(0,T;H_0^1(\Omega))\cap H^1(0,T;L^2(\Omega)),\\
&u_{N^\prime}\!-\!a_{0,N^\prime}\!-\!t a_{1,N^\prime} \rightharpoonup \widetilde u 
&&\quad \mbox{weakly in}\ H_\alpha(0,T;H^{-1}(\Omega)),
\end{aligned}
\right.
\end{equation}
as $N^\prime\to \infty$. 
Since $a_{0,N^\prime}\rightarrow a_0$ and $a_{1,N^\prime}\rightarrow a_1$ strongly in 
$L^2(\Omega)$ as $N^\prime\to \infty$, we have 
$u_{N^\prime}-a_{0,N^\prime}-t a_{1,N\prime} \rightharpoonup u-a_0-t a_1$
weakly in $L^2(0,T;H_0^1(\Omega))\cap H^1(0,T;L^2(\Omega))$. 
By noting that the weak limit is unique, we obtain $\widetilde u = u-a_0-t a_1$. 
Therefore, we derive
\begin{align*}
&\quad \|u\|_{L^2(0,T;H_0^1(\Omega))} + \|u\|_{H^1(0,T;L^2(\Omega))} 
+ \|u-a_0-t a_1\|_{H_\alpha(0,T;H^{-1}(\Omega))} \\
&\le \liminf_{N^\prime\to\infty} \bigg(\|u_{N^\prime}\|_{L^2(0,T;H_0^1(\Omega))} 
+ \|u_{N^\prime}\|_{H^1(0,T;L^2(\Omega))} \\
&\quad + \|u_{N^\prime}-a_{0,N^\prime}-t a_{1,N\prime}\|_{H_\alpha(0,T;H^{-1}(\Omega))} \bigg)\\
&\le C \left(\|a_0\|_{H_0^1(\Omega)} + \|a_1\|_{L^2(\Omega)} + \|F\|_{L^2(0,T;L^2(\Omega))}\right).
\end{align*}
In terms of \eqref{eq:ap7}, noting that $L^{\infty}(0,T;H^1_0(\OOO))
= (L^1(0,T;H^{-1}(\OOO)))^*$: the adjoint space, we see that 
$\sup_{N\in \N} \Vert u_N\Vert_{(L^1(0,T;H^{-1}(\OOO)))^*} < \infty$.
Therefore, by the Banach-Alaoglu-Bourbaki theorem (e.g., Chapter 3 in
Brezis \cite{Br}), we can extract a subsequence of $\N$, denoted by the same 
letter $N'$, such that there exists $\widetilde{u} 
\in (L^1(0,T;H^{-1}(\OOO)))^*$ and $u_{N'} \longrightarrow 
\widetilde{u}$ weakly$^*$ converges in $(L^1(0,T;H^{-1}(\OOO)))^*$, that is,
$$
_{L^\infty(0,T;H^1_0(\OOO)))}{\langle u_{N'},\, \psi
\rangle}_{L^1(0,T;H^{-1}(\OOO))}
\, \longrightarrow \,\,
_{L^\infty(0,T;H^1_0(\OOO)))}{\langle \widetilde{u},\, \psi
\rangle}_{L^1(0,T;H^{-1}(\OOO))},
$$
for any $\psi \in L^1(0,T;H^{-1}(\OOO))$.  By the uniqueness of the limit,
we obtain $\widetilde{u} = u$ and 
$$
\Vert u\Vert_{(L^1(0,T;H^{-1}(\OOO)))^*}
\le \liminf_{N' \to\infty} \Vert u_{N'} \Vert_{(L^1(0,T;H^1_0(\OOO)))^*},
$$
see e.g., Chapter 3 in \cite{Br}.
Again applying \eqref{eq:ap7}, we see
$$
\Vert u\Vert_{L^\infty(0,T;H^1_0(\OOO)))}
\le C(\Vert a_0\Vert_{H^1_0(\OOO)} + \Vert a_1\Vert_{L^2(\OOO)}
+ \Vert F\Vert_{L^2(0,T;L^2(\OOO))}).
$$
In a similar way, we obtain the estimate of $\partial_t^{\alpha-1}(u-a_0)$ as follows:
$$
\|\partial_t^{\alpha-1}(u-a_0)\|_{L^\infty(0,T;L^2(\Omega))} 
\le C \left(\|a_0\|_{H_0^1(\Omega)} + \|a_1\|_{L^2(\Omega)} + \|F\|_{L^2(0,T;L^2(\Omega))}\right).
$$

Next, we show that $u$ is a weak solution to \eqref{equ:pde2}. 
Recall that $u$ is a weak solution to \eqref{equ:pde2} if 
$u-a_0-t a_1\in H_\alpha(0,T;H^{-1}(\Omega))$, $u\in L^2(0,T;H_0^1(\Omega))$ and 
\begin{equation}
\label{eq:def-sol}
_{H^{-1}(\Omega)}\langle \partial_t^\alpha (u-a_0-ta_1), \psi \rangle
_{H_0^1(\Omega)} 
= B[u(t),\psi;t] + {_{H^{-1}(\Omega)}\langle F(\cdot,t), \psi \rangle
_{H_0^1(\Omega)} },
\end{equation}
for almost all $t\in [0,T]$ and any $\psi\in H_0^1(\Omega)$. 
Here $B[u,v;t]$ is a bilinear form defined by 
$$
B[u,v;t] := \int_\Omega \Big( -\sum_{i,j=1}^n a_{ij}(t)(\partial_i u)
\partial_j v 
+ \sum_{j=1}^n b_j(t)(\partial_j u)\,v + c(t)uv \Big) dx.
$$
Recall that $u_N$ solves \eqref{equ:pde2b} in the $L^2$-sense. 
Then for any $\widetilde \psi\in L^2(0,T;H_0^1(\Omega))$, we multiply \eqref{equ:pde2b} by $\widetilde \psi$ 
and then integrate over $\Omega\times (0,T)$, which gives
\begin{align*}
&\quad \int_0^T {_{H^{-1}(\Omega)}\langle \partial_t^\alpha (u_N-a_{0,N}-ta_{1,N}), \widetilde \psi(t) 
\rangle_{H_0^1(\Omega)} } dt\\
&= \int_0^T B[u_N(t),\widetilde\psi(t);t] dt 
+ \int_0^T {_{H^{-1}(\Omega)}\langle F_N(\cdot,t), \widetilde\psi(t) \rangle_{H_0^1(\Omega)} } dt,
\end{align*}
where we used integration by parts: 
\begin{align*}
\int_0^T\! (-Au_N(t), \widetilde \psi(t))_{L^2(\Omega)} dt 
&= \!\int_0^T\!\! \int_\Omega \!
\Big(\sum_{i,j=1}^n\! \partial_j (a_{ij}\partial_i u_N) \!+\! \sum_{j=1}^n b_j \partial_j u_N \!+\! cu_N\Big)\widetilde \psi 
dxdt\\
&= \int_0^T B[u_N(t),\widetilde \psi(t);t] dt.
\end{align*}
By taking the subsequence $\{u_{N^\prime}\}$ satisfying 
\eqref{eq:convergence}, we have
\begin{align}
\nonumber
&\quad \int_0^T {_{H^{-1}(\Omega)}\langle \partial_t^\alpha (u_{N^\prime}-a_{0,N^\prime}-ta_{1,N^\prime}), 
\widetilde \psi(t) \rangle_{H_0^1(\Omega)} } dt\\
\label{eq:exist1}
&= \int_0^T B[u_{N^\prime}(t),\widetilde \psi(t);t] dt 
+ \int_0^T {_{H^{-1}(\Omega)}\langle F_{N^\prime}(\cdot,t), \widetilde \psi(t) \rangle_{H_0^1(\Omega)} } dt.
\end{align}
According to the weak convergence \eqref{eq:convergence}, we have 
$u_{N^\prime} \rightharpoonup u$ weakly in $L^2(0,T;H_0^1(\Omega))$ and 
$\partial_t^\alpha (u_{N^\prime}-a_{0,N^\prime}-t a_{1,N^\prime}) \rightharpoonup 
\partial_t^\alpha (u-a_0-t a_1)$ weakly in $L^2(0,T;H^{-1}(\Omega))$ as $N^\prime\to\infty$. 
Moreover,  
$F_{N^\prime} \rightarrow F$ strongly in $L^2(0,T;L^2(\Omega))$ which implies 
$F_{N^\prime} \rightarrow F$ also weakly in $L^2(0,T;H^{-1}(\Omega))$. 
Since $\int_0^T B[\cdot,\widetilde \psi;t]dt$ is a bounded linear functional on $L^2(0,T;H_0^1(\Omega))$, 
we let $N^\prime\to\infty$ in \eqref{eq:exist1} and obtain 
\begin{align*}
&\quad \int_0^T {_{H^{-1}(\Omega)}\langle \partial_t^\alpha (u-a_0-ta_1), 
\widetilde \psi(t) \rangle_{H_0^1(\Omega)} } dt\\
&= \int_0^T B[u(t),\widetilde \psi(t);t] dt 
+ \int_0^T {_{H^{-1}(\Omega)}\langle F(\cdot,t), \widetilde \psi(t) \rangle_{H_0^1(\Omega)} } dt.
\end{align*}
Note that in the above equality, $\widetilde \psi$ can be taken arbitrarily in $L^2(0,T;H_0^1(\Omega))$. 
Therefore, \eqref{eq:def-sol} holds true for any $\psi\in H_0^1(\Omega)$ and 
almost all $t\in [0,T]$, and 
thus, $u$ is a weak solution to \eqref{equ:pde2}. 

\noindent {\bf Step 4. Uniqueness of solution}

In the case of $\alpha\le 1$, as one can see from the proof of Theorem 4.1 in \cite{KRY20}, 
the uniqueness of solution can be easily proved by taking $\psi=u(t)\in H_0^1(\Omega)$ 
in \eqref{eq:def-sol} and applying the generalized Gr\"onwall's inequality. 
However, in our case of $1<\alpha\le 2$, we cannot do in the same way 
because $u$ has not enough regularity as the test function. 
Here we modify the idea used in \cite[Chapter 7]{E98} to show the uniqueness of the solution 
by proving first the uniqueness of the integral of the solution.   

By taking the difference of two possible solutions, it is sufficient to show that the weak solution 
$\overline u\in L^2(0,T;H_0^1(\Omega))$ to 
$$
\left\{
\begin{aligned}
& \partial_t^\alpha \overline u(x,t) = -A(x,t) \overline u(x,t), &&\quad (x,t)\in \Omega\times (0,T),\\
& \overline u \in H_\alpha(0,T;H^{-1}(\Omega)), 
\end{aligned}
\right.
$$
must be zero. 
By the definition of the weak solution, we readily obtain 
$$
\int_0^\tau {_{H^{-1}(\Omega)}\langle \partial_s^\alpha \overline u(s), \psi \rangle_{H_0^1(\Omega)} } ds
= \int_0^\tau B[\overline u(s),\psi;s] ds,
$$
for any $\psi\in H_0^1(\Omega)$ and almost all $\tau\in [0,T]$. 
Now we take $\psi = \overline u(\tau)\in H_0^1(\Omega)$, and for arbitrarily fixed $t\in [0,T]$, 
we integrate the above equality over $\tau\in (0,t)$, so that 
\begin{align}
\label{eq:uni1}
\int_0^t \int_0^\tau {_{H^{-1}(\Omega)}\langle \partial_s^\alpha \overline u(s), \overline u(\tau) 
\rangle_{H_0^1(\Omega)} } ds d\tau
= \int_0^t \int_0^\tau B[\overline u(s),\overline u(\tau);s] ds d\tau,
\end{align}
for any $t\in [0,T]$. 
Next, we estimate both sides of \eqref{eq:uni1}. 

By using
\begin{align*}
\int^{\tau}_0
{_{H^{-1}(\Omega)}\langle \partial_s^\alpha \overline u(s), \overline u(\tau) 
\rangle_{H_0^1(\Omega)} } ds
&= {_{H^{-1}(\Omega)}\left\langle \int^{\tau}_0
\partial_s^\alpha \overline u(s) ds, \, \, \overline u(\tau) 
\right\rangle_{H_0^1(\Omega)} } \\
&= {_{H^{-1}(\Omega)}\langle (J^1\partial_{\tau}^\alpha \overline u)(\tau), \,
\overline u(\tau) \rangle_{H_0^1(\Omega)} },
\end{align*}
and employing Theorem \ref{thm:J4}(iii), we have
\begin{align*}
&\int_0^t \int_0^\tau {_{H^{-1}(\Omega)}\langle \partial_s^\alpha \overline u(s), 
\overline u(\tau) \rangle_{H_0^1(\Omega)} } ds d\tau 
= \int_0^t {_{H^{-1}(\Omega)}\langle J^1\partial_\tau^\alpha \overline u(\tau), 
\overline u(\tau) \rangle_{H_0^1(\Omega)} } d\tau  \\
&= \int_0^t {_{H^{-1}(\Omega)}\langle \partial_\tau^{\alpha-1} \overline u(\tau), 
\overline u(\tau) \rangle_{H_0^1(\Omega)} } d\tau
= (J^1({_{H^{-1}(\Omega)}\langle \partial_\tau^{\alpha-1} \overline u, 
\overline u \rangle_{H_0^1(\Omega)} }) )(t)\\
& = (J^{2-\alpha}(J^{\alpha-1} ({_{H^{-1}(\Omega)}\langle 
\partial_\tau^{\alpha-1} \overline u, \overline u 
\rangle_{H_0^1(\Omega)} })))(t)\\
&\ge \frac12 J^{2-\alpha}\|\overline u(\cdot,t)\|_{L^2(\Omega)}^2. 
\end{align*}
For the final inequality, we used Lemma \ref{lem:coer} in Appendix \ref{subsec:app3} 
and the fact that $J^{2-\alpha}$ preserves the ordering, 
that is, $J^{2-\alpha}w(t) \ge 0$ for $0\le t \le T$ if 
$w(t) \ge 0$ for $0\le t \le T$.

By letting $\overline v = J^1 \overline u$, which implies $\overline u = \frac{\partial \overline v}{\partial t}$,
the above inequality yields
\begin{equation}
\label{eq:uni2}
\int_0^t \int_0^\tau {_{H^{-1}(\Omega)}\langle \partial_s^\alpha \overline u(s), \overline u(\tau) 
\rangle_{H_0^1(\Omega)} } ds d\tau 
\ge \frac12 J^{2-\alpha} \left\|\frac{\partial \overline v}{\partial t} (\cdot,t)\right\|_{L^2(\Omega)}^2. 
\end{equation}
On the other hand, we let
\begin{align*}
\int_0^t \!\int_0^\tau\! B[\overline u(s),\overline u(\tau);s] ds d\tau 
&= - \int_0^t \!\int_0^\tau \!\int_\Omega \sum_{i,j=1}^n a_{ij}(s) 
(\partial_i \overline u(s)) \partial_j \overline u(\tau) 
dx ds d\tau\\
&\quad+ \!\int_0^t \!\int_0^\tau \!\int_\Omega 
\!\Big(\sum_{j=1}^n b_j(s)\partial_j \overline u(s) + c(s)\overline u(s)\Big)
\overline u(\tau) dx ds d\tau\\
&=: I_1(t) + I_2(t).
\end{align*}
By noting that $\overline v(0) = J^1 \overline u(0) =0$, we use Fubini's theorem and integration by parts 
to estimate
\begin{align*}
I_1(t) &= - \int_\Omega \sum_{i,j=1}^n \int_0^t a_{ij}(s) \partial_i \overline u(s) 
\int_s^t \partial_j \overline u(\tau) d\tau dsdx \\
&= - \int_\Omega \sum_{i,j=1}^n \int_0^t a_{ij}(s) (\partial_i \overline u(s)) 
\partial_j \overline v(t) ds dx 
\!+\! \int_\Omega \sum_{i,j=1}^n \int_0^t a_{ij}(s) (\partial_i \overline u(s)) 
\partial_j \overline v(s) ds dx\\
&= - \int_\Omega \sum_{i,j=1}^n a_{ij}(t) (\partial_i \overline v(t)) 
\partial_j \overline v(t) dx 
+ \int_\Omega \sum_{i,j=1}^n \int_0^t \frac{\ppp a_{ij}}{\ppp s}(s) 
(\partial_i \overline v(s)) \partial_j \overline v(t) ds dx\\
&\quad + \frac12 \int_\Omega \sum_{i,j=1}^n a_{ij}(t) (\partial_i \overline v(t))
 \partial_j \overline v(t) dx
\!-\! \frac12 \int_\Omega \sum_{i,j=1}^n \int_0^t \frac{\ppp a_{ij}}{\ppp s}
(s) 
(\partial_i \overline v(s))\partial_j \overline v(s)ds dx\\
&= -\frac12 \int_\Omega \sum_{i,j=1}^n a_{ij}(t) (\partial_i \overline v(t))
 \partial_j \overline v(t) dx
+ \int_\Omega \sum_{i,j=1}^n \int_0^t \frac{\ppp a_{ij}}{\ppp s}(s) 
(\partial_i \overline v(s)) \partial_j \overline v(t) ds dx\\
&\quad - \frac12 \int_\Omega \sum_{i,j=1}^n \int_0^t \frac{\ppp a_{ij}}{\ppp s}
(s) (\partial_i \overline v(s)) 
\partial_j \overline v(s) ds dx.
\end{align*}
Then by the assumption (3) on $a_{ij}$ and H\"older's inequality, we have
\begin{align*}
&-\frac12 \int_\Omega \sum_{i,j=1}^n a_{ij}(t) (\partial_i \overline v(t))
 \partial_j \overline v(t) dx
\le -\frac{\sigma_0}{2}\|\nabla \overline v(\cdot,t)\|_{L^2(\Omega)}^2,\\
&\int_\Omega \!\sum_{i,j=1}^n \!\int_0^t \! \frac{\ppp a_{ij}}{\ppp s}(s) 
(\partial_i \overline v(s)) \partial_j \overline v(t) ds dx
\le \varepsilon \|\nabla \overline v(\cdot,t)\|_{L^2(\Omega)}^2 
\!+\! \frac{C}{\varepsilon}\int_0^t \|\nabla \overline v(\cdot,s)\|_{L^2(\Omega)}^2 ds,\\
&- \frac12 \int_\Omega \sum_{i,j=1}^n \int_0^t \frac{\ppp a_{ij}}{\ppp s}(s) 
(\partial_i \overline v(s)) \partial_j \overline v(s) ds dx
\le C\int_0^t \|\nabla \overline v(\cdot,s)\|_{L^2(\Omega)}^2 ds.
\end{align*}
Thus, we obtain
$$
I_1(t) \le -\left(\frac{\sigma_0}{2} - \varepsilon\right) 
\|\nabla \overline v(\cdot,t)\|_{L^2(\Omega)}^2 
+ \left(C + \frac{C}{\varepsilon}\right)\int_0^t \|\nabla \overline v(\cdot,s)\|
_{L^2(\Omega)}^2 ds.
$$
Similarly, we estimate
\begin{align*}
I_2(t) &= \int_\Omega \int_0^t \Big(\sum_{j=1}^n b_j(s)\partial_j \overline u(s) + c(s)\overline u(s)\Big) 
\int_s^t \overline u(\tau) d\tau dsdx \\
&= \int_\Omega \int_0^t \Big(\sum_{j=1}^n b_j(s)\partial_j \overline u(s) + c(s)\overline u(s)\Big) \overline v(t) ds dx\\ 
&\quad - \int_\Omega \int_0^t \Big(\sum_{j=1}^n b_j(s)\partial_j \overline u(s) + c(s)\overline u(s)\Big) \overline v(s) ds dx \\
&= - \int_\Omega \int_0^t \Big(\sum_{j=1}^n \frac{\ppp b_j}{\ppp s}(s)
  \partial_j \overline v(s) 
+ \frac{\ppp c}{\ppp s}(s) \overline v(s)\Big) \overline v(t) ds dx \\
&\quad + \int_\Omega \int_0^t \Big(\sum_{j=1}^n \frac{\ppp b_j}{\ppp s}(s) 
 \partial_j \overline v(s) 
+ \frac{\ppp c}{\ppp s}(s) \overline v(s)\Big) \overline v(s) ds dx \\
&\quad + \int_\Omega \int_0^t \Big(\sum_{j=1}^n b_j(s)\partial_j \overline v(s) + c(s)\overline v(s)\Big) 
\frac{\partial \overline v}{\partial s}(s) ds dx.
\end{align*}
By the assumptions on $b_j,c$, H\"older's inequality and the Poincar\'e inequality, we have 
$$
I_2(t) \le \varepsilon \|\nabla \overline v(\cdot,t)\|_{L^2(\Omega)}^2 
+ \left(C \!+\! \frac{C}{\varepsilon}\right)\!\int_0^t\! \|\nabla \overline 
v(\cdot,s)\|_{L^2(\Omega)}^2 ds 
+ C\!\int_0^t\! \left\| \frac{\partial \overline v}{\partial s} (\cdot,s)\right\|_{L^2(\Omega)}^2 
ds. 
$$
Hence by taking $\varepsilon= \frac{\sigma_0}{8}$, we obtain
\begin{align}
\nonumber
&\quad \int_0^t \int_0^\tau B[\overline u(s),\overline u(\tau);s] ds d\tau \\
\label{eq:uni3}
&\le -\frac{\sigma_0}{4}\|\nabla \overline v(\cdot,t)\|_{L^2(\Omega)}^2
+ C\int_0^t \|\nabla \overline v(\cdot,s)\|_{L^2(\Omega)}^2 ds 
+ C\int_0^t \left\|\frac{\partial \overline v}{\partial s}(\cdot,s)\right\|_{L^2(\Omega)}^2 ds.
\end{align}
Inserting \eqref{eq:uni2} and \eqref{eq:uni3} into \eqref{eq:uni1} yields
\begin{align*}
&\quad J^{2-\alpha}\left\| \frac{\partial \overline v}{\partial t}(\cdot,t)
\right\|_{L^2(\Omega)}^2 + \|\nabla \overline v(\cdot,t)\|_{L^2(\Omega)}^2\\
&\le C\int_0^t \left\| \frac{\partial \overline v}{\partial s}(\cdot,s)\right\|
_{L^2(\Omega)}^2 ds
+ C\int_0^t \|\nabla \overline v(\cdot,s)\|_{L^2(\Omega)}^2 ds.
\end{align*}
In the same way as we did in Step 3, we have
$$
\int_0^t \left\| \frac{\partial \overline v}{\partial s}(\cdot,s)\right\|
_{L^2(\Omega)}^2 ds 
= J^{\alpha-1}J^{2-\alpha}\left\| \frac{\partial \overline v}{\partial t}(\cdot,t)
\right\|_{L^2(\Omega)}^2,
$$
and
$$
\int_0^t \|\nabla \overline v(\cdot,s)\|_{L^2(\Omega)}^2 ds
\le \Gamma(\alpha-1)T^{2-\alpha} J^{\alpha-1}\|\nabla \overline v(\cdot,t)\|
_{L^2(\Omega)}^2, 
$$
and thus, with a new constant $C>0$, we obtain
\begin{align*}
&\quad J^{2-\alpha}\left\| \frac{\partial \overline v}{\partial t}(\cdot,t)\right\|
_{L^2(\Omega)}^2 
+ \|\nabla \overline v(\cdot,t)\|_{L^2(\Omega)}^2\\
&\le CJ^{\alpha-1}\left(J^{2-\alpha}\left\|\frac{\partial \overline v}{\partial t}(\cdot,t)
\right\|_{L^2(\Omega)}^2 
+ \|\nabla \overline v(\cdot,t)\|_{L^2(\Omega)}^2\right).
\end{align*}
Therefore, we apply the generalized Gr\"onwall's inequality ($\gamma=\alpha-1$ in Lemma \ref{lem:Gronwall} in Section 5) and 
the Poincar\'e inequality to reach $\overline v(\cdot,t) = 0$ for almost all
$t\in [0,T]$. 
Noting that $\overline u = \frac{\partial}{\partial t} \overline v$, we have also $\overline u(\cdot,t) = 0$ 
for almost all $t\in [0,T]$. 
\qed

\section{Proof of Theorem \ref{thm:pde2}}
\label{sec:pde2}
In this section, we prove Theorem \ref{thm:pde2}. 
If $a_1\in H^2(\Omega)\cap H_0^1(\Omega)$ and we allow $\|a_1\|_{H^2(\Omega)}$ 
on the right-hand side of the regularity estimate, then the proof would be easy 
since one could denote $v=u-a_0-t a_1$ and apply $L^2$ estimate to $v$ in space.  
In Theorem \ref{thm:pde2}, we have only $\|a_1\|_{H_0^1(\Omega)}$ on the right-hand side, 
and thus we need some technical treatments as follows. 

If we prove the a priori estimate 
\begin{align*}
&\quad \|\partial_t^\alpha (u_N-a_{0,N}-t a_{1,N})\|_{L^\infty(0,T;L^2(\Omega))} 
+ \|\partial_t^{\alpha-1}(u_N-a_{0,N})\|_{L^\infty(0,T;H_0^1(\Omega))} \\
&\quad + \|u_N\|_{H^1(0,T;H_0^1(\Omega))} 
+ \|u_N\|_{L^\infty(0,T;H^2(\Omega))} \\
&\le C\left(\|a_0\|_{H^2(\Omega)} + \|a_1\|_{H_0^1(\Omega)} + \|F\|_{H^1(0,T;L^2(\Omega))}\right),
\end{align*}
then we can derive the desired regularity and estimate for $u$ by a similar argument used in Step 3 of 
the proof in Section \ref{sec:pde}. Then it is sufficient to establish such a uniform estimate. 

Here we consider slightly different approximate solutions due to some technical reasons. 
For simplicity, we still use the same notation $u_N$. 
Precisely, we construct the approximate solutions:
$$
u_N(x,t) = \sum_{k=1}^N p_k^N(t) \varphi_k(x), \quad (x,t)\in \Omega\times (0,T),
$$
satisfying
\begin{equation}
\label{equ:pde3b}
\left\{
\begin{aligned}
& \partial_t^\alpha (u_N - a_{0,N} - t a_{1,N}) 
= -A(x,t) (u_N-a_{0,N})(x,t) + G_N(x,t), \\
& u_N - a_{0,N} - t a_{1,N} \in H_\alpha(0,T;L^2(\Omega)), 
\end{aligned}
\right.
\end{equation}
where
\begin{align*}
&a_{j,N} = \sum_{k=1}^N a_k^j \varphi_k, \quad 
a_k^j := (a_j,\varphi_k)_{L^2(\Omega)}, \quad j=0,1,\\
&G_N(t) = \sum_{k=1}^N g_k(t)\varphi_k, \quad 
g_k(t) := (F(t)-A(t)a_0,\varphi_k)_{L^2(\Omega)}, \quad 0<t<T.
\end{align*}
Then $p^N = (p_1^N,\ldots,p_N^N)^T$ solves 
\begin{equation}
\label{equ:pde3c}
\left\{
\begin{aligned}
& \partial_t^\alpha (p^N \!-\! a^0 \!-\! t a^1) = Q(t)(p^N(t) \!-\! a^0) + g(t) = Q(t)p^N(t) + (g(t)\!-\!Q(t)a^0), \\
& p^N - a^0 - t a^1 \in (H_\alpha(0,T))^N,
\end{aligned}
\right.
\end{equation}
where $Q(t)=(q_{\ell k}(t))_{k,\ell=1}^N$ is defined by \eqref{def:Q} and $g=(g_1,\ldots,g_N)^T$. 
By the assumption $F-Aa_0\in H_1(0,T;L^2(\Omega))$, we find $g\in (H_1(0,T))^N$, and thus,
$(g-Qa^0) + Qa^0 = g\in (H_1(0,T))^N$. 
Hence we can employ Corollary \ref{coro:ode4} and obtain 
$p^N - a^0 - t a^1 \in (H_{\alpha+1}(0,T))^N$, which implies 
$$
u_N - a_{0,N} - t a_{1,N} \in H_{\alpha+1}(0,T;H^2(\Omega)\cap H_0^1(\Omega)). 
$$
Moreover, since $t a_{1,N}\in H_1(0,T; H^2(\Omega)\cap H_0^1(\Omega))$, immediately we have 
$u_N - a_{0,N}\in H_1(0,T;H^2(\Omega)\cap H_0^1(\Omega))$. 
Then we take the time derivative in \eqref{equ:pde3b} and we obtain
$$
\frac{\partial}{\partial t}\partial_t^\alpha (u_N - a_{0,N} - t a_{1,N}) 
= -A \frac{\partial}{\partial t}(u_N-a_{0,N}) + \frac{\partial G_N}{\partial t} - \frac{\partial A}{\partial t} (u_N-a_{0,N}).
$$
Now by setting
$v := \frac{\partial}{\partial t} (u_N-a_{0,N}-t a_{1,N}) \in H_\alpha(0,T;H^2(\Omega)\cap H_0^1(\Omega))$, 
we rewrite the above equation by
\begin{equation}
\label{eq:im1}
\partial_t^\alpha v(t) = -Av(t) -A(t)a_{1,N} + \frac{\partial G_N}{\partial t}(t) - \frac{\partial A}{\partial t} (u_N-a_{0,N}),
\quad 0<t<T.  
\end{equation}
We multiply \eqref{eq:im1} by $\partial_t^{\alpha-1} v$ and 
integrate over $\Omega$. 
By integration by parts and substituting $t$ by $s$, we obtain
\begin{align}
\nonumber
&\quad \frac12 \frac{d}{d s} \|\partial_s^{\alpha-1} v(\cdot,s)\|_{L^2(\Omega)}^2 
+ \int_\Omega \sum_{i,j=1}^n a_{ij}(x,s) (\partial_j v(x,s)) 
\partial_s^{\alpha-1} \partial_i v(x,s) dx\\
\nonumber
&= \!\int_\Omega \!\Big(\sum_{j=1}^n b_j(s) \partial_j v(s) \!+\! c(s)v(s)\Big) \partial_s^{\alpha-1} v(s) dx
\!-\! \int_\Omega \!\sum_{i,j=1}^n a_{ij}(s)(\partial_j a_{1,N})
 \partial_s^{\alpha-1} \partial_i v(s) dx\\
\nonumber
&\quad + \int_\Omega \Big(\sum_{j=1}^n b_j(s) \partial_j a_{1,N} + c(s) a_{1,N}\Big) 
\partial_s^{\alpha-1} v(s) dx
+ \int_\Omega \frac{\partial G_N}{\partial s}(s) \partial_s^{\alpha-1} v(s) dx\\
\label{eq:im2}
&\quad - \int_\Omega \frac{\partial A}{\partial s}(s)
 (u_N(s)-a_{0,N}) \partial_s^{\alpha-1} v(s) dx
=: \sum_{k=1}^5 I_k(s).
\end{align}
After the second line of \eqref{eq:im2}, we remark that we omit the variable
$x$ for concise descriptions, and $v(t)$ means $v(x,t)$ or $v(\cdot,t)$. 
  
We estimate the right-hand side of \eqref{eq:im2} term by term. 
Recalling that the coefficients $a_{ij},b_j,c$ are bounded from above, 
by H\"older's inequality and the Poincar\'e inequality, we have
\begin{align*}
&I_1(s) \le \varepsilon\|\nabla v(\cdot,s)\|_{L^2(\Omega)}^2 
+ \frac{C}{\varepsilon}\|\partial_s^{\alpha-1} v(\cdot,s)\|_{L^2(\Omega)}^2\quad 
\mbox{for any }\varepsilon>0,\\
&I_3(s) \le C\left(\|a_{1,N}\|_{H_0^1(\Omega)}^2 
+ \|\partial_s^{\alpha-1} v(\cdot,s)\|_{L^2(\Omega)}^2\right),\\
&I_4(s) \le C\left(\left\|\frac{\partial G_N}{\partial s}(\cdot,s)\right\|_{L^2(\Omega)}^2 
+ \|\partial_s^{\alpha-1} v(\cdot,s)\|_{L^2(\Omega)}^2\right),\\
&I_5(s) \le C\left(\|u_N(s)-a_{0,N}\|_{H^2(\Omega)}^2 
+ \|\partial_s^{\alpha-1} v(\cdot,s)\|_{L^2(\Omega)}^2\right).
\end{align*}
Moreover, by the equation \eqref{equ:pde3b} and the elliptic regularity 
(e.g., \cite[Theorem 4, Chapter VI]{E98}), we obtain
\begin{align}
\nonumber
&\quad \|u_N(s)-a_{0,N}\|_{H^2(\Omega)}^2 \\
\nonumber
&\le C\|-A(s)(u_N(s)-a_{0,N})\|_{L^2(\Omega)}^2 + C\|u_N(s)-a_{0,N}\|_{L^2(\Omega)}^2\\
\nonumber
&\le C\|\partial_s^\alpha (u_N(s)\!-\!a_{0,N}\!-\!s a_{1,N})\|_{L^2(\Omega)}^2 
+ C\|G_N(s)\|_{L^2(\Omega)}^2 + C\|u_N(s)\!-\!a_{0,N}\|_{L^2(\Omega)}^2\\
\label{eq:im3}
&\le C\|\partial_s^{\alpha-1} v(\cdot,s)\|_{L^2(\Omega)}^2 + C\|G_N(s)\|_{L^2(\Omega)}^2
+ C\|u_N(s)-a_{0,N}\|_{L^2(\Omega)}^2, 
\end{align}
and hence
$$
I_5(s) \le C\left(\|\partial_s^{\alpha-1} v(\cdot,s)\|_{L^2(\Omega)}^2 + \|G_N(s)\|_{L^2(\Omega)}^2
+ \|u_N(s)-a_{0,N}\|_{L^2(\Omega)}^2\right).
$$
Furthermore, by Proposition \ref{prop:expression}, 
we estimate $I_2$ as follows:
\begin{align*}
&I_2(s) = - \int_\Omega \sum_{i,j=1}^n a_{ij}(s)(\partial_j a_{1,N}) 
\frac{\partial}{\partial s} J^{2-\alpha} \partial_i v(s) dx\\
&= -\frac{d}{ds} \!\int_\Omega \sum_{i,j=1}^n \!a_{ij}(s)
(\partial_j a_{1,N}) J^{2-\alpha} \partial_i v(s) dx
\!+\! \int_\Omega \sum_{i,j=1}^n \!\frac{\partial a_{ij}}{\partial s}(s)
(\partial_j a_{1,N}) J^{2-\alpha} \partial_i v(s) dx\\
&\le -\frac{d}{ds} \!\int_\Omega \sum_{i,j=1}^n \!a_{ij}(s)
(\partial_j a_{1,N}) J^{2-\alpha} \partial_i v(s) dx 
\!+\! CJ^{2-\alpha} \left(\|a_{1,N}\|_{H_0^1(\Omega)}^2 \!+\! \|\nabla v\|_{L^2(\Omega)}^2\right)\\
&\le -\frac{d}{ds} \!\int_\Omega \sum_{i,j=1}^n \!a_{ij}(s)(\partial_j a_{1,N})
 J^{2-\alpha} \partial_i v(s) dx
\!+\! CT^{2-\alpha}\|a_{1,N}\|_{H_0^1(\Omega)}^2 \!+\! CJ^{2-\alpha}\|\nabla v\|_{L^2(\Omega)}^2.
\end{align*}
Then we apply the operator $J^1$ to \eqref{eq:im2}, that is, 
we integrate \eqref{eq:im2} over $s\in (0,t)$. 
Noting that $\partial_s^{\alpha-1}v(0)=0$ and $J^1 = J^{2-\alpha} J^{\alpha-1}$, we derive
\begin{align*}
&\quad \frac12 \|\partial_t^{\alpha-1} v(\cdot,t)\|_{L^2(\Omega)}^2 
+ J^{2-\alpha} J^{\alpha-1}
\int_\Omega \sum_{i,j=1}^n a_{ij}(s) \partial_j v(s) 
\partial_s^{\alpha-1} \partial_i v(s) dx\\
&\le -\int_\Omega \sum_{i,j=1}^n a_{ij}(t)\partial_j a_{1,N} J^{2-\alpha}\partial_i v(t) dx \\
&\quad + C\left(\|a_{1,N}\|_{H_0^1(\Omega)}^2 
+ \|G_N\|_{L^2(0,t;L^2(\Omega))}^2 
+ \left\|\frac{\partial G_N}{\partial t}\right\|_{L^2(0,t;L^2(\Omega))}^2 \right)\\
&\quad + \varepsilon J^1\|\nabla v\|_{L^2(\Omega)}^2 + CJ^{3-\alpha}\|\nabla v\|_{L^2(\Omega)}^2 \\
&\quad + \left(C+\frac{C}{\varepsilon}\right)J^1\|\partial_s^{\alpha-1} v\|_{L^2(\Omega)}^2 
+ CJ^1\|u_N-a_{0,N}\|_{L^2(\Omega)}^2.
\end{align*}
By H\"older's inequality and the order preserving property of 
$J^{2-\alpha}$, for any $\varepsilon>0$, we have
\begin{align*}
-\int_\Omega \sum_{i,j=1}^n a_{ij}(t)(\partial_j a_{1,N})
 J^{2-\alpha}\partial_i v(t) dx
&\le CJ^{2-\alpha} \int_\Omega \sum_{i,j=1}^n |\partial_j a_{1,N}| |\partial_i v| dx\\
&\le \frac{C}{\varepsilon}\|a_{1,N}\|_{H_0^1(\Omega)}^2 
+ \varepsilon J^{2-\alpha}\|\nabla v\|_{L^2(\Omega)}^2.
\end{align*}
Also we employ the coercivity inequality (Lemma \ref{lem:coer2} in Appendix) to obtain
\begin{align*}
&\quad J^{2-\alpha} J^{\alpha-1}\int_\Omega \!\sum_{i,j=1}^n a_{ij}(s) 
(\partial_j v(s)) \partial_s^{\alpha-1} \partial_i v(s) dx \\
&\ge \frac{\sigma_0}{2}J^{2-\alpha} \|\nabla v\|_{L^2(\Omega)}^2 
- CJ^{3-\alpha} \|\nabla v\|_{L^2(\Omega)}^2. 
\end{align*}
Moreover, we have
\begin{align*}
J^1\|\nabla v\|_{L^2(\Omega)}^2 &= \int_0^t \|\nabla v(\cdot,s)\|_{L^2(\Omega)}^2 ds \\
&= \int_0^t \Gamma(3-\alpha)(t-s)^{\alpha-2} \frac{(t-s)^{2-\alpha}}{\Gamma(3-\alpha)} 
\|\nabla v(\cdot,s)\|_{L^2(\Omega)}^2 ds \\
&\le \Gamma(3-\alpha)T^{\alpha-2} J^{3-\alpha}\|\nabla v\|_{L^2(\Omega)}^2.
\end{align*}
Note that the above inequality is trivial if $\alpha=2$.
Therefore, by taking $\varepsilon>0$ small enough, the above four estimates yield
\begin{align}
\nonumber
&\quad \|\partial_t^{\alpha-1} v(\cdot,t)\|_{L^2(\Omega)}^2 
+ J^{2-\alpha}\|\nabla v\|_{L^2(\Omega)}^2\\
\nonumber
&\le C\!\left(\|a_{1,N}\|_{H_0^1(\Omega)}^2 \!+\! \|G_N\|_{L^2(0,t;L^2(\Omega))}^2 
\!+\! \left\|\frac{\partial G_N}{\partial t}\right\|_{L^2(0,t;L^2(\Omega))}^2 \!+\! J^1\|u_N-a_{0,N}\|_{L^2(\Omega)}^2\right)\\
\label{eq:im4}
&\quad + CJ^1\left(\|\partial_s^{\alpha-1} v\|_{L^2(\Omega)}^2 
+ J^{2-\alpha}\|\nabla v\|_{L^2(\Omega)}^2\right),
\end{align}
for $0\le t\le T$. 
It is readily to verify that 
\begin{align*}
&\|a_{1,N}\|_{H_0^1(\Omega)}^2 \le \|a_1\|_{H_0^1(\Omega)}^2,\\
&\|G_N\|_{L^2(0,t;L^2(\Omega))}^2 
\le C\left(\|F\|_{L^2(0,T;L^2(\Omega))}^2 + \|a_0\|_{H^2(\Omega)}^2\right),\\
&\left\|\frac{\partial G_N}{\partial t}\right\|_{L^2(0,t;L^2(\Omega))}^2 
\le C\left(\|F\|_{H^1(0,T;L^2(\Omega))}^2 + \|a_0\|_{H^2(\Omega)}^2\right).
\end{align*}
According to the estimate of Theorem \ref{thm:pde1}, we have
\begin{align*}
J^1\|u_N-a_{0,N}\|_{L^2(\Omega)}^2 
&\le C\|u_N\|_{L^2(0,t;L^2(\Omega))}^2 + C\|a_{0,N}\|_{L^2(\Omega)}^2\\
&\le C\left(\|F\|_{L^2(0,T;L^2(\Omega))}^2 + \|a_0\|_{H_0^1(\Omega)}^2 
+ \|a_1\|_{L^2(\Omega)}^2\right),
\end{align*}
in $(0,T)$. Combining the above estimates with \eqref{eq:im4}, we reach
\begin{align*}
&\quad \|\partial_t^{\alpha-1} v(\cdot,t)\|_{L^2(\Omega)}^2 
+ J^{2-\alpha}\|\nabla v\|_{L^2(\Omega)}^2\\
&\le C\left(\|a_0\|_{H^2(\Omega)}^2 + \|a_1\|_{H_0^1(\Omega)}^2 
+ \|F\|_{H^1(0,T;L^2(\Omega))}^2\right) \\
&\quad + CJ^1\left(\|\partial_s^{\alpha-1} v\|_{L^2(\Omega)}^2 
+ J^{2-\alpha}\|\nabla v\|_{L^2(\Omega)}^2\right),
\end{align*}
for $0\le t\le T$. 
By Gr\"onwall's inequality ($\gamma=1$ in 
Lemma \ref{lem:Gronwall} in Section 5), 
the above inequality implies 
\begin{align}
\nonumber
&\quad \|\partial_t^{\alpha-1} v(\cdot,t)\|_{L^2(\Omega)}^2 
+ J^{2-\alpha}\|\nabla v\|_{L^2(\Omega)}^2 \\
\label{eq:im5}
&\le C\left(\|a_0\|_{H^2(\Omega)}^2 + \|a_1\|_{H_0^1(\Omega)}^2 
+ \|F\|_{H^1(0,T;L^2(\Omega))}^2\right)
\end{align}
for $0\le t\le T$. Then we return to \eqref{eq:im3} and see that
\begin{align*}
\|u_N(t)\|_{H^2(\Omega)}^2 
&\le C\|u_N(t)-a_{0,N}\|_{H^2(\Omega)}^2 + C\|a_{0,N}\|_{H^2(\Omega)}^2 \\
&\le C\left(\|a_0\|_{H^2(\Omega)}^2 + \|a_1\|_{H_0^1(\Omega)}^2 
+ \|F\|_{H^1(0,T;L^2(\Omega))}^2\right). 
\end{align*}
Finally, substituting $v=\frac{\partial}{\partial t} (u_N-u_{0,N}-t a_{1,N})$ into \eqref{eq:im5} yields
\begin{align*}
&\quad \|\partial_t^{\alpha} (u_N(t)\!-\!u_{0,N}\!-\!t a_{1,N})\|_{L^2(\Omega)}^2 
\!+\! J^{2-\alpha}\left\|\nabla 
\left(\frac{\partial}{\partial s} u_N-a_{1,N}\right)
\right\|_{L^2(\Omega)}^2 
\!+\! \|u_N(t)\|_{H^2(\Omega)}^2\\
&\le C\left(\|a_0\|_{H^2(\Omega)}^2 + \|a_1\|_{H_0^1(\Omega)}^2 
+ \|F\|_{H^1(0,T;L^2(\Omega))}^2\right).
\end{align*}
We apply the operator $J^{\alpha-1}$ to the second term on the left-hand side and take $t=T$,  
which implies
$$
\left\|\nabla \frac{\partial}{\partial t} u_N\right\|_{L^2(0,T;L^2(\Omega))}^2 
\le C\left(\|a_0\|_{H^2(\Omega)}^2 + \|a_1\|_{H_0^1(\Omega)}^2 
+ \|F\|_{H^1(0,T;L^2(\Omega))}^2\right).
$$
In a similar way as we derive \eqref{eq:ap4}, we obtain
\begin{align}
\nonumber
&\|\nabla \partial_t^{\alpha-1} (u_N-a_{0,N}-ta_{1,N})\|_{L^\infty(0,T;L^2(\Omega))}^2 \\
\label{eq:rev1}
&\le C\left(\|a_0\|_{H^2(\Omega)}^2 + \|a_1\|_{H_0^1(\Omega)}^2 
+ \|F\|_{H^1(0,T;L^2(\Omega))}^2\right).
\end{align}
Since we can use the triangle inequality to estimate
\begin{align*}
&\quad \|\nabla \partial_t^{\alpha-1} (u_N-a_{0,N})\|_{L^\infty(0,T;L^2(\Omega))}^2 \\
&\le 2\|\nabla \partial_t^{\alpha-1} (u_N-a_{0,N}-ta_{1,N})\|_{L^\infty(0,T;L^2(\Omega))}^2 + \frac{2T^{2(2-\alpha)}}{(\Gamma(3-\alpha))^2}\|\nabla a_{1,N}\|_{L^2(\Omega)}^2,
\end{align*}
this and \eqref{eq:rev1} imply
\begin{align*}
\|\nabla \partial_t^{\alpha-1}\! (u_N-a_{0,N})\|_{L^\infty(0,T;L^2(\Omega))}^2 
\!\le C\left(\|a_0\|_{H^2(\Omega)}^2 \!+\! \|a_1\|_{H_0^1(\Omega)}^2 
\!+\! \|F\|_{H^1(0,T;L^2(\Omega))}^2\right).
\end{align*}
This completes the proof of Theorem \ref{thm:pde2}. 
\qed
\section{Appendix}

\subsection{Characterization of $H_{\gamma}(0,T;\R)$ in terms of 
Sobolev-Slobodeckij spaces}
\label{subsec:app1}

Although we do not directly use it, we provide a characterization 
of $H_{\gamma}(0,T):= H_{\gamma}(0,T;\R)$ 
in terms of Sobolev-Slobodeckij spaces \cite{S58}.

For $0<\gamma<1$, 
we define the Sobolev-Slobodeckij space 
$H^{\gamma}(0,T) := 
\{ u \in L^2(0,T);\, \Vert u\Vert_{H^{\gamma}(0,T)} 
< \infty\}$ of real-valued functions as follows.
$$
\|u\|_{H^\gamma(0,T)} := 
\left(\|u\|_{L^2(0,T)}^2 
+ \int_0^T \int_0^T \frac{\vert u(t)-u(s)\vert^2}
{|t-s|^{1+2\gamma}} dt ds\right)^{\frac{1}{2}}.
$$
For $\gamma>1$ and $\gamma\not\in \N$, we write $\gamma = \ell 
+ \theta$ where $\ell \in \N$ and $\theta\in (0,1)$, we define 
$$
H^{\gamma}(0,T) := \left\{ u\in H^{\ell}(0,T);\, 
\frac{d^{\ell}}{dt^{\ell}} u
\in H^{\theta}(0,T) \right\},
$$ 
and it is known that this is a Banach space with respect to the norm 
$$
\|u\|_{H^\gamma(0,T)} := 
\left(\|u\|_{H^{\ell}(0,T)}^2 
+ \left\| \frac{d^{\ell}}{dt^{\ell}} u\right\|
_{H^{\theta}(0,T)}^2 \right)^{\frac{1}{2}}.
$$
We set $H^0(0,T):= L^2(0,T)$ and $\frac{d^0u}{dt^0} := u$.
Then
\begin{theorem}
\label{thm:char}
For $0<\gamma\le 1$, we have 
$$
H_\gamma(0,T) =
\left\{
\begin{aligned}
&\{u\in H^\gamma(0,T);\, u(0) = 0\} \quad \mbox{if 
$\frac{1}{2}<\gamma\le 1$}, \\
&\left\{u\in H^{\frac{1}{2}}(0,T);\, \int_0^T t^{-1}\vert u(t)\vert^2 dt 
< \infty\right\} \quad \mbox{if $\gamma= \frac{1}{2}$}, \\
& H^\gamma(0,T) \quad \mbox{if $0<\gamma<\frac{1}{2}$},
\end{aligned}
\right.
$$
and the norm $\Vert u\Vert_{H_{\gamma}(0,T)} 
:= \Vert J^{-\gamma}u\Vert_{L^2(0,T)}$ 
is equivalent to  
$$
\left\{
\begin{aligned}
& \|u\|_{H^\gamma(0,T)}, \quad 0<\gamma\le 1, \gamma\not= \frac{1}{2},\\
& \left(\|u\|_{H^{\frac{1}{2}}(0,T)}^2 
+ \int_0^T t^{-1}\vert u(t)\vert^2 dt\right)^{\frac{1}{2}}, 
\quad \gamma=\frac{1}{2}.
\end{aligned}
\right.
$$
Moreover, for $\gamma = \ell+\theta$ with $\ell\in \mathbb{N}$ and $0<\theta<1$, we have
$$
H_{\ell+\theta}(0,T) = \left\{u\in H_\ell(0,T); \, \frac{d^\ell}{dt^\ell} 
u\in H_\theta(0,T)\right\},
$$
and the norm $\|u\|_{H_{\ell+\theta}(0,T)}$ is equivalent to 
$$
\left(\|u\|_{H^\ell(0,T)}^2 + \left\|\frac{d^\ell}{dt^\ell} u\right\|
_{H_\theta(0,T)}^2\right)^{\frac{1}{2}}.
$$
\end{theorem}

This characterization indicates that $H_\gamma(0,T;\mathbb{R})$ 
is a subspace of 
the corresponding Sobolev-Slobodeckij space $H^\gamma(0,T;\mathbb{R})$. 

The proof of the theorem in Kubica, Ryszewska and Yamamoto \cite{KRY20} 
is based on 
\begin{enumerate}
\item[(i)]
We estimate $\Vert \ppp_t^1u\Vert_{L^2(0,T)}$ by 
$\Vert S^{\frac{1}{2}}u\Vert_{L^2(0,T)}$.  
Here the operator $S$ is defined by 
$$
(Su)(t) = - \frac{d^2u}{dt^2}(t), \quad 0<t<T, \quad
\mathcal{D}(S) = \left\{u \in H^2(0,T);\, u(0) = \frac{du}{dt}(T) = 0\right\},
$$
and $S^{\frac{1}{2}}$ is the square root of $S$. 

\item[(ii)]
We show that $\ppp_t^{\gamma}$ with the domain $H_{\gamma}(0,T)$ is the 
fractional power (Tanabe \cite{Ta}) of the operator $\ppp_t^1$ with the domain 
$H_1(0,T)$.

\item[(iii)]
Applying the Heinz-Kato inequality (e.g., \cite{Ta}), we obtain
$$
\mathcal{D}((S^{\frac{1}{2}})^{\gamma}) = \mathcal{D}(\ppp_t^{\gamma})
= H_{\gamma}(0,T).
$$

\item[(iv)]
Applying the interpolation inequality (\cite{LM}) based on the spectral
representation of $S^{\frac{1}{2}}$, we can characterize 
$\mathcal{D}((S^{\frac{1}{2}})^{\gamma})$ by means of
the Sobolev-Slobodeckij spaces.
\end{enumerate}
Thus, for general Banach space $X$, we can expect to prove the 
corresponding result to Theorem 8, but we here omit the details.

We can refer to Yamamoto \cite{Yam1,Yam2} as for 
different approaches about the characterization of $H_{\gamma}(0,T) =
\mathcal{D}(\ppp_t^{\gamma})$.

\subsection{Proof of Lemma \ref{lem:dense}}
\label{subsec:app2}

{\bf Part (i)}
Let $0 \le t_1 < t \le T$.
We have
\begin{align*}
 \Gamma(\alpha)(J^{\alpha}u(t_1) - J^{\alpha}u(t))
&= - \int^t_{t_1} (t-s)^{\alpha-1} u(s) ds\\
&\quad + \int^{t_1}_0 ((t_1-s)^{\alpha-1} - (t-s)^{\alpha-1}) u(s) ds
=: I_1 + I_2.
\end{align*}
We estimate $I_1$ by the Cauchy-Schwarz inequality and $\alpha > \frac{1}{2}$:
\begin{align*}
\Vert I_1\Vert_X &\le \int^t_{t_1} (t-s)^{\alpha-1} \Vert u(s)\Vert_X ds\\
&\le \left( \int^t_{t_1} (t-s)^{2\alpha-2} ds\right)^{\frac{1}{2}}
   \left( \int^t_{t_1} \Vert u(s)\Vert_X^2 ds\right)^{\frac{1}{2}}
\le \left( \frac{(t-t_1)^{2\alpha-1}}{2\alpha-1}\right)^{\frac{1}{2}}
\Vert u\Vert_{\LLLLLL}.
\end{align*}
As for $I_2$ we argue as follows. By $\alpha > \frac{1}{2}$, we can find 
small $\ep > 0$ such that $\alpha > \frac{1}{2} + \ep$.
Since $\vert a^{1-\alpha} - b^{1-\alpha} \vert 
\le \vert a-b\vert^{1-\alpha}$ for all $a, b \ge 0$ by $1-\alpha < 1$, we 
see
$$
\vert (t-s)^{1-\alpha} - (t_1-s)^{1-\alpha} \vert 
\le \vert (t-s) - (t_1-s) \vert^{1-\alpha}
= \vert t-t_1\vert^{1-\alpha}.
$$
Therefore, using $t-s \ge t-t_1$ and $t-s \ge t_1-s$ for 
$0 < s < t_1 < t$, we obtain
\begin{align*}
\Vert I_2\Vert_X &= \left\Vert
\int^{t_1}_0 \vert (t_1-s)^{\alpha-1} - (t-s)^{\alpha-1} \vert 
u(s) ds \right\Vert_X \\
&\le \int^{t_1}_0 \vert (t_1-s)^{\alpha-1} - (t-s)^{\alpha-1}\vert 
\Vert u(s)\Vert_X ds \\
&= \int^{t_1}_0 
\frac{ \vert (t-s)^{1-\alpha} - (t_1-s)^{1-\alpha}\vert}
{(t-s)^{1-\alpha} (t_1-s)^{1-\alpha}} \Vert u(s)\Vert_X ds\\
&\le \int^{t_1}_0 \frac{
\vert (t-t_1)^{1-\alpha}}
{(t-s)^{1-\alpha-\ep}(t-s)^{\ep} (t_1-s)^{1-\alpha}} \Vert u(s)\Vert_X ds \\
&\le \frac{(t-t_1)^{1-\alpha}}{(t-t_1)^{1-\alpha-\ep}}
\int^{t_1}_0 (t_1-s)^{\alpha-1-\ep} \Vert u(s)\Vert_X ds \\
&\le (t-t_1)^{\ep}\left( \int^{t_1}_0 (t_1-s)^{2\alpha-2-2\ep} ds 
\right)^{\frac{1}{2}}
\left( \int^{t_1}_0 \Vert u(s)\Vert_X ds\right)^{\frac{1}{2}} \\
&\le (t-t_1)^{\ep} \left( \frac{t_1^{2\alpha-1-2\ep}}{2\alpha-1-2\ep}
\right)^{\frac{1}{2}} \Vert u\Vert_{\LLLLLL}.
\end{align*}
Here we used $\alpha > \frac{1}{2} + \ep$.
Hence, 
$$
\Vert \Gamma(\alpha)(J^{\alpha}u(t) - J^{\alpha}u(t_1))\Vert_X
\le C((t-t_1)^{\alpha-\frac{1}{2}} + (t-t_1)^{\ep}T^{\alpha-\frac{1}{2}-\ep})
\Vert u\Vert_{\LLLLLL},
$$
which means that $J^{\alpha}u \in C([0,T];X)$.

Next, by $\alpha>\frac{1}{2}$, we apply the Cauchy-Schwarz inequality to have
\begin{align*}
\Vert \Gamma(\alpha)J^{\alpha}u(t)\Vert_X
&= \left\Vert \int^t_0 (t-s)^{\alpha-1} u(s) ds \right\Vert_X
\le \left( \int^t_0 (t-s)^{2\alpha-2} ds \right)^{\frac{1}{2}}
\left( \int^t_0 \Vert u(s)\Vert_X ds \right)^{\frac{1}{2}}\\
&\le \left( \frac{t^{2\alpha-1}}{2\alpha-1}\right)^{\frac{1}{2}}
\Vert u\Vert_{\LLLLLL}
\le Ct^{2\alpha-1} \Vert u\Vert_{\LLLLLL}.
\end{align*}
Therefore, $\lim_{t\downarrow 0} \Vert J^{\alpha}u(t)\Vert_X = 0$.
Thus, the proof of (i) is complete.

\noindent {\bf Part (ii)}
Let $\alpha > \frac{3}{2}$.  We can write $\alpha = 1 + \beta$ where
$\beta > \frac{1}{2}$.  Let $u \in H_{\alpha}(0,T;X)$.
Then $u = J^{1+\beta}w$ with some $w \in \LLLLLL$.
By Lemma \ref{lem:J1}(iii), we have $u = J^1J^{\beta}w$, which implies
$\frac{d}{dt}u = \frac{d}{dt} J^1J^{\beta}w = J^{\beta}w$. 
In terms of $\beta > \frac{1}{2}$, part (i) yields $\frac{d}{dt}u \in C([0,T];X)$ 
and $\frac{d}{dt}u(0)=0$. 
Thus, the proof of (ii) is complete.

\noindent {\bf Part (iii)} 
We set $\, _{0}C^1([0,T];X) := \{ u\in C^1([0,T];X);\, u(0) = 0\}$.
Let $u \in \Halp$ and $\ep>0$ be given arbitrarily.  Then we can 
find $w \in \LLLLLL$ such that $u = J^{\alpha} w$.
Since $C^{\infty}_0(0,T;X)$ is dense in $\LLLLLL$ by the mollifier
(e.g., \cite{AF75}), there exists $\varphi_{\ep} \in C^{\infty}_0(0,T;X)$ such 
that
\begin{equation}
\label{app:eq1}
\Vert \varphi_{\ep} - w \Vert_{\LLLLLL} < \ep.  
\end{equation}
Noting that 
$$
(J^{\alpha}\varphi_{\ep})(t) = \frac{1}{\Gamma(\alpha)} \int^t_0 (t-s)^{\alpha-1}
\varphi_{\ep}(s) ds
= \frac{1}{\Gamma(\alpha)} \int^t_0 s^{\alpha-1}\varphi_{\ep}(t-s) ds,
$$
we can readily verify that $J^{\alpha}\varphi_{\ep} \in \, _{0}C^1([0,T];X)$.

Now, by the definition of $\Vert \cdot\Vert_{\Halp}$, we have
$$
\Vert u - J^{\alpha}\varphi_{\ep}\Vert_{\Halp} 
= \Vert J^{\alpha}w - J^{\alpha}\varphi_{\ep}\Vert_{\Halp}
= \Vert w - \varphi_{\ep} - w\Vert_{\LLLLLL}.
$$
Therefore, \eqref{app:eq1} implies $\Vert u - J^{\alpha}\varphi_{\ep}\Vert_{\Halp}
< \ep$.  Since $J^{\alpha}\varphi_{\ep} \in \, _{0}C^1([0,T];X)$,
this means that $u \in \ooo{\, _{0}C^1([0,T];X)}^{\Halp}$.
Thus, the proof of (iii) is complete.
\qed

\subsection{Coercivity inequalities of the time-fractional derivative}
\label{subsec:app3}

Here we show two coercivity inequalities, which are 
used several times in the proofs of the main results. 
\begin{lemma}[First coercivity inequality]
\label{lem:coer}
Let $0<\gamma\le 1$. For each $u\in L^2(0,T;H_0^1(\Omega))\cap H_\gamma(0,T;H^{-1}(\Omega))$, 
\begin{align*}
J^{\gamma}({_{H^{-1}(\Omega)}\langle\partial_t^\gamma u(s),u(s)\rangle
_{H_0^1(\Omega)}})(t) 
&= \int_0^t \frac{(t-s)^{\gamma-1}}{\Gamma(\gamma)}
{_{H^{-1}(\Omega)}\langle\partial_t^\gamma u(s),u(s)\rangle_{H_0^1(\Omega)}} ds \\
&\ge \frac12 \|u(\cdot,t)\|_{L^2(\Omega)}^2.
\end{align*}
\end{lemma}
For this lemma, in terms of Lemma \ref{lem:dense}(iii), it is sufficient to 
prove for $u\in {_{0}C^1([0,T];L^2(\OOO))}$, and then 
we can follow the proof of Theorem 3.4(\romannumeral2) of \cite{KRY20}. 
Here we omit the details.

\begin{lemma}[Second coercivity inequality]
\label{lem:coer2}
Let $0<\gamma\le 1$ and $a_{ij}=a_{ji}\in W^{1,\infty}(\Omega\times (0,T))$, $i,j=1,\ldots,n$ 
satisfy the assumption \eqref{assump:aij}. Then for each 
$v=(v_1,\ldots,v_n)^T\in H_1(0,T;L^2(\Omega))$,
\begin{align*}
&\quad \int_0^t \frac{(t-s)^{\gamma-1}}{\Gamma(\gamma)} 
\int_\Omega \sum_{i,j=1}^n a_{ij}(x,s) v_i(x,s) \partial_s^\gamma v_j(x,s) dx ds\\
&\ge \frac{\sigma_0}{2} \|v(\cdot,t)\|_{L^2(\Omega)}^2 - C\int_0^t \|v(\cdot,s)\|_{L^2(\Omega)}^2 ds,
\end{align*}
for some constant $C=C(\gamma,\|a_{ij}\|_{W^{1,\infty}})>0$. 
\end{lemma}
\begin{proof}
For completeness, we give the detailed proof but we also refer to \cite[Proof of Theorem 4.2]{KRY20} 
for a similar idea. For the case of $\gamma=1$, the left-hand side of our desired inequality reads
$$
\int_0^t \int_\Omega \sum_{i,j=1}^n a_{ij}(x,s) v_i(x,s) 
\frac{\ppp v_j}{\ppp s}(x,s) dx ds.
$$
Then by the assumption \eqref{assump:aij} and integration by parts, 
we can immediately prove the desired inequality. 
Now we consider the case of $0<\gamma<1$. Let
$$
I_0(x,s) := \sum_{i,j=1}^n a_{ij}(x,s) v_i(x,s) \partial_s^\gamma v_j(x,s), 
\quad (x,s)\in \Omega\times (0,T). 
$$
In the following estimate of $I_0$, we may omit $x$ for simplicity. 
Since we assume $v\in H_1(0,T;L^2(\Omega))$, 
the fractional derivative $\partial_s^\gamma$ coincides with the Caputo fractional derivative 
$_C\partial_s^\gamma$ and by $a_{ij}=a_{ji}$, $i,j=1,\ldots,n$, we have
\begin{align*}
I_0(s) &= \sum_{i,j=1}^n a_{ij}(s) v_i(s) 
\int_0^s \frac{(s-\tau)^{-\gamma}}{\Gamma(1-\gamma)} \frac{\partial}{\partial \tau} v_j(\tau) d\tau\\
&= \sum_{i,j=1}^n a_{ij}(s) \int_0^s \frac{(s-\tau)^{-\gamma}}
{\Gamma(1-\gamma)} 
\left( \frac{\partial v_j(\tau)}{\ppp \tau}\right) v_i(\tau) d\tau \\
&\quad + \sum_{i,j=1}^n a_{ij}(s) \int_0^s \frac{(s-\tau)^{-\gamma}}{\Gamma(1-\gamma)} 
\left( \frac{\partial v_j(\tau)}{\ppp\tau}\right) (v_i(s)-v_i(\tau)) d\tau\\
&=: I_1(s) + I_2(s).
\end{align*}
By noting $v_i(0)=0$ for $i=1,\ldots,n$ and
\begin{align*}
&\quad \lim_{\tau\to s} \frac{(s-\tau)^{-\gamma}}{\Gamma(1-\gamma)} (v_i(s)-v_i(\tau)) (v_j(s)-v_j(\tau))\\
&\le \lim_{\tau\to s} \frac{(s-\tau)^{1-\gamma}}{\Gamma(1-\gamma)} \|v_i\|_{H^1(0,T)} \|v_j\|_{H^1(0,T)}
=0,
\end{align*}
we use integration by parts to estimate $I_2$ as follows:
\begin{align*}
I_2(s) &= - \frac12 \sum_{i,j=1}^n a_{ij}(s) \int_0^s \frac{(s-\tau)^{-\gamma}}{\Gamma(1-\gamma)} 
\frac{\partial}{\partial \tau} \left((v_j(s)-v_j(\tau))(v_i(s)-v_i(\tau))\right) d\tau\\
&= \frac12 \sum_{i,j=1}^n a_{ij}(s) \frac{s^{-\gamma}}{\Gamma(1-\gamma)} 
\left(v_j(s)v_i(s)\right) \\
&\quad + \frac{\gamma}{2} \sum_{i,j=1}^n a_{ij}(s) \int_0^s \frac{(s-\tau)^{-\gamma-1}}{\Gamma(1-\gamma)} 
\left((v_j(s)-v_j(\tau))(v_i(s)-v_i(\tau))\right) d\tau\\
&\ge \frac{\sigma_0}{2} \frac{s^{-\gamma}}{\Gamma(1-\gamma)} 
|v(s)|_{\R^n}^2
+ \frac{\sigma_0}{2} \frac{\gamma}{\Gamma(1-\gamma)} 
\int_0^s \frac{|v(s)-v(\tau)|_{\R^n}^2}{(s-\tau)^{\gamma+1}} d\tau
\ge 0.
\end{align*}
Here in the last line we used the assumption \eqref{assump:aij}. 
Similarly, we estimate $I_1$ as follows:
\begin{align*}
I_1(s) &= \frac12 \sum_{i,j=1}^n a_{ij}(s) \int_0^s \frac{(s-\tau)^{-\gamma}}{\Gamma(1-\gamma)} 
\frac{\partial}{\partial \tau} (v_i(\tau) v_j(\tau)) d\tau \\
&= \frac12 \int_0^s \frac{(s-\tau)^{-\gamma}}{\Gamma(1-\gamma)} 
\frac{\partial}{\partial \tau} \Big(\sum_{i,j=1}^n (a_{ij}(s)-a_{ij}(\tau)) v_i(\tau) v_j(\tau)\Big) d\tau \\
&\quad + \frac12 \int_0^s \frac{(s-\tau)^{-\gamma}}{\Gamma(1-\gamma)} 
\frac{\partial}{\partial \tau} \Big(\sum_{i,j=1}^n a_{ij}(\tau) v_i(\tau) v_j(\tau)\Big) d\tau\\
&= -\frac{\gamma}{2} 
\int_0^s\! \sum_{i,j=1}^n \frac{(s-\tau)^\gamma}{\Gamma(1-\gamma)} 
\frac{a_{ij}(s)\!-\!a_{ij}(\tau)}{s-\tau} v_i(\tau) v_j(\tau) d\tau 
\!+\! \frac12 J^{1-\gamma} \frac{\partial}{\partial s} \Big(\sum_{i,j=1}^n a_{ij} v_i v_j\Big).
\end{align*}
Since $a_{ij}\in W^{1,\infty}(\Omega\times (0,T))$, there exists a generic constant $C>0$ 
depending on some norm of $a_{ij}$ such that 
$$
\left|\frac{a_{ij}(s)-a_{ij}(\tau)}{s-\tau}\right| \le C.
$$
Thus, by the triangle inequality $|AB|\le \frac12(A^2+B^2)$, we have the lower bound of $I_1$:
$$
I_1(s) \ge \frac12 J^{1-\gamma} \frac{\partial}{\partial s} \Big(\sum_{i,j=1}^n a_{ij} v_i v_j\Big) 
- C J^{1-\gamma} |v|_{\mathbb{R}^n}^2(s).
$$
Finally, by Fubini's theorem, we obtain
\begin{align*}
\int_0^t \frac{(t-s)^{\gamma-1}}{\Gamma(\gamma)} \int_\Omega I_0(x,s) dxds
&=J^{\gamma}\int_\Omega I_0 dx \\
&\ge \frac12 \int_\Omega \sum_{i,j=1}^n a_{ij}(t) v_i(t) v_j(t) dx 
- C \int_0^t \int_\Omega |v|_{\mathbb{R}^n}^2 dxds\\
&\ge \frac{\sigma_0}{2} \|v(\cdot,t)\|_{L^2(\Omega)}^2 - C\int_0^t \|v(\cdot,s)\|_{L^2(\Omega)}^2 ds.
\end{align*}
\end{proof}

\subsection{Generalized Gr\"onwall's inequality}
\label{subsec:app4}
Here, we state the following generalized Gr\"onwall's inequality. 
\begin{lemma}
\label{lem:Gronwall}
Let $0<\gamma\le 1$ and $r\in L^1(0,T)$ fulfill $r\ge 0$ in $(0,T)$. 
We assume that $u\in L^1(0,T)$ satisfies
$$
0\le u(t) \le r(t) + C\int_0^t (t-s)^{\gamma-1} u(s) ds, \quad 0\le t\le T,
$$
for some positive constant $C>0$. Then there exist positive constants $C_1,C_2>0$, 
depending on $\gamma$ and $C$ but are independent on $T$, such that
$$
u(t) \le r(t) + C_1e^{C_2t}\int_0^t(t-s)^{\gamma-1} r(s) ds, \quad 0\le t\le T. 
$$
Moreover, if $r$ is a continuous non-decreasing function, then there exists a constant $C_3>0$, 
depending on $C$, $\gamma$ and $T$, such that
$$
u(t) \le C_3 r(t), \quad 0\le t\le T.
$$ 
\end{lemma}
The proof for $0<\gamma<1$ can be found in e.g., \cite[Appendix A]{KRY20}, 
while the case of $\gamma=1$ is the well-known Gr\"onwall's inequality. 

\section*{Acknowledgements}
The authors sincerely thank the anonymous referee for the valuable comments to help improve this article, in particular, Section 1.2. 
The first author was supported by Grant-in-Aid for 
JSPS (Japan Society for the Promotion of Science) Fellows 20F20319. 
The second author was supported by 
Grant-in-Aid for Scientific Research (A) 20H00117 and Grant-in-Aid 
for Challenging Research (Pioneering) 21K18142 of Japan Society for the Promotion of Science.


\begin{thebibliography}{}
%
%

\bibitem{AF75}
Adams, R.A., Fournier, J.J.F.: 
{\em Sobolev Spaces}. 
Vol. 140, 2nd edition, Academic Press, San Diego (2003).

\bibitem{Br}
Brezis, H.:
{\em 
Functional Analysis, Sobolev Spaces and Partial Differential Equations}.
Springer, New York (2011).

\bibitem{E96}
Engler, H.: 
Global smooth solutions for a class of parabolic integrodifferential equations. 
Trans. Amer. Math. Soc. \textbf{348}, 267--290 (1996). 
https://doi.org/10.1090/S0002-9947-96-01472-9

\bibitem{E98}
Evans, L.: 
{\em Partial Differential Equations}. second edition
American Mathematical Society, Providence, Rhode Island (2010).


\bibitem{GKMR14}
Gorenflo, R., Kilbas, A.A., Mainardi, F., Rogosin, S.V.: 
{\em Mittag-Leffler Functions, Related Topics and Applications}. 
Springer Monographs in Mathematics, Springer, Berlin (2014). 
https://doi.org/10.1007/978-3-662-43930-2


\bibitem{GLS90}
Gripenberg, G., Londen, S.O., Staffans, O.: 
{\em Volterra Integral and Functional Equations}. 
Cambridge University Press (1990). 
https://doi.org/10.1017/CBO9780511662805

\bibitem{GLY15}
Gorenflo, R., Luchko, Y., Yamamoto, M.:  
Time-fractional diffusion equation in the fractional Sobolev spaces. 
Fract. Calc. Appl. Anal. \textbf{18}, 799--820 (2015). 
https://doi.org/10.1515/fca-2015-0048


\bibitem{HKP20}
Han, B.-S., Kim, K.-H., Park, D.:
Weighted $L_q(L_p)$-estimate with Muckenhoupt weights for the diffusion-wave equations with time-fractional derivatives. 
J. Differ. Equ. \textbf{269}, 3515--3550 (2020).
https://doi.org/10.1016/j.jde.2020.03.005.

\bibitem{KY}
Kian, K., Yamamoto,M.:
Well-posedness for weak and strong solutions of non-homogeneous 
initial boundary value problems for fractional diffusion equations.
Fract. Calc. Appl. Anal. {\bf 24}, 168-201 (2021).
https://doi.org/10.1515/fca-2021-0008

\bibitem{K08}
Kochubei, A. N.: 
Distributed order calculus and equations of ultraslow diffusion.   
J. Math. Anal. Appl. \textbf{340}, 252--281 (2008). 
https://doi.org/10.1016/j.jmaa.2007.08.024

\bibitem{KRY20}
Kubica, A., Ryszewska, K., Yamamoto, M.:  
{\em Time-Fractional Differential Equations: A Theoretical Introduction}.  
Springer, Singapore (2020).
https://doi.org/10.1007/978-981-15-9066-5

\bibitem{LLY15}
Li, Z., Liu, Y., Yamamoto, M.: 
Initial-boundary value problems for multi-term time-fractional diffusion equations with positive 
constant coefficients. 
Applied Mathematics and Computation \textbf{257}, 381--397 (2015).
https://doi.org/10.1016/j.amc.2014.11.073

\bibitem{LM}
Lions, J.-L., Magenes, E.: 
{\em Non-homogeneous Boundary Value Problems and Applications}. Vol. 2. 
Springer, Berlin (1972). 
https://doi.org/10.1007/978-3-642-65217-2

\bibitem{L11}
Luchko, Y.: 
Initial-boundary-value problems for the generalized multi-term time-fractional diffusion equation.
J. Math. Anal. Appl. \textbf{374}, 538--548 (2011). 
https://doi.org/10.1016/j.jmaa.2010.08.048

\bibitem{LY16}
Luchko, Y., Yamamoto, M.:  
General time-fractional diffusion equation: some uniqueness and existence results for the 
initial-boundary-value problems.  
Fract. Calc. Appl. Anal. \textbf{19}, 676--695 (2016).
https://doi.org/10.1515/fca-2016-0036

\bibitem{M95}
Mainardi, F.: 
Fractional diffusive waves in viscoelastic solids.  
In: Wegner, J. L., Norwood, F. R. (eds.) Nonlinear Waves in Solids, pp. 93--97. 
ASME/AMR, Fairfield (1995).

\bibitem{MMP08}
Mainardi, F., Mura, A., Pagnini, G., Gorenflo, R.:  
Time-fractional diffusion of distributed order. 
Journal of Vibration and Control \textbf{14}, 1267--1290 (2008). 
https://doi.org/10.1177/1077546307087452

\bibitem{N86}
Nigmatullin, R.R.:  
The realization of the generalized transfer equation in a medium with fractal geometry.  
Physica Status Solidi \textbf{133}, 425--430 (1986). 
https://doi.org/10.1002/pssb.2221330150

\bibitem{P23}
Park, D.: 
Weighted maximal $L_q(L_p)$-regularity theory for time-fractional diffusion-wave equations with variable coefficients. 
J. Evol. Equ. \textbf{23}, 12 (2023). 
https://doi.org/10.1007/s00028-022-00866-8

\bibitem{P99}
Podlubny, I.: 
{\em Fractional Differential Equations}. 
Academic Press, San Diego (1999).

\bibitem{P93}
Pr\"uss, J.: 
{\em Evolutionary Integral Equations and Applications}.
Birkh\"auser, Basel (1993).
https://doi.org/10.1007/978-3-0348-8570-6

\bibitem{SY11}
Sakamoto, K., Yamamoto, M.:  
Initial value/boundary value problems for fractional diffusion-wave equations and applications to 
some inverse problems. 
J. Math. Anal. Appl. \textbf{382}, 426--447 (2011).
https://doi.org/10.1016/j.jmaa.2011.04.058

\bibitem{SKM93}
Samko, S.G., Kilbas, A.A., Marichev, O.I.:
{\em Fractional Integrals and Derivatives: Theory and Applications}. 
Gordon and Breach, Amsterdam (1993).

\bibitem{S58}
Slobodeckij, L.N.:  
Generalized Sobolev spaces and their applications to boundary value problems of partial 
differential equations.  
Leningrad. Gos. Ped. Inst. Ucep. Zap. \textbf{197}, 54--112 (1958).

\bibitem{Ta}
Tanabe, H.: {\em Equations of Evolution}. Pitman, London (1979).

\bibitem{Yagi}
Yagi, A.: 
{\em Abstract Parabolic Evolution Equations and their Applications}. 
Springer Berlin (2010).
https://doi.org/10.1007/978-3-642-04631-5

\bibitem{Yam1}
Yamamoto, M.:
Weak solutions to non-homogeneous boundary value problems for
time-fractional diffusion equations.
J. Math. Anal. Appl. \textbf{460}, 365--381 (2018).

\bibitem{Yam2}
Yamamoto, M.:
{\em Fractional derivatives and time-fractional ordinary differential equations in $L^p$-space}. 
Preprint (2022), arXiv:2201.07094

\bibitem{Z06}
Zacher, R.:
Quasilinear parabolic integro-differential equations with nonlinear boundary conditions. 
Differential Integral Equations \textbf{19}, 1129--1156 (2006). 
https://doi.org/10.57262/die/1356050312

\bibitem{Z09}
Zacher, R.:  
Weak solutions of abstract evolutionary integro-differential equations in Hilbert spaces.  
Funkcial. Ekvac. \textbf{52}, 1--18 (2009).
https://doi.org/10.1619/fesi.52.1

\end{thebibliography}


\end{document}